\documentclass{amsart}

\usepackage{amssymb}
\usepackage{amsbsy}
\usepackage{amscd}
\usepackage{amsmath}
\usepackage{amsthm}
\usepackage{amsfonts,mathrsfs}

\newtheorem{theorem}{Theorem}
\newtheorem{lemma}[theorem]{Lemma}
\newtheorem{proposition}{Proposition}

\newtheorem{defn}{Definition}\numberwithin{defn}{section}
\newtheorem{example}{Example}

\newcommand{\N}{\mathbb{N}}
\newcommand{\R}{\mathbb{R}}

\newcommand{\E}{\mathbb{E}}
\newcommand{\PP}{\mathbb{P}}

\newcommand{\cH}{\mathcal{H}}
\newcommand{\cM}{\mathcal M}

\newcommand{\cU}{\mathscr{U}}

\newcommand{\F}{\mathcal{F}}
\newcommand{\B}{\mathcal{B}}
\newcommand{\LL}{\mathcal{L}}
\newcommand{\NN}{\mathcal{N}}
\newcommand{\EE}{\mathcal{E}}
\newcommand{\FF}{\mathscr{F}}

\newcommand{\gga}{\gamma}            

\newcommand{\LLog}{L\log^{1/2}\!\!L(X,\gamma)}
\newcommand{\aplim}{\mathop{\rm ap\, lim}}
\newcommand{\one}{{\bf 1}}

\newcommand{\FCb}{{\mathscr F}C^1_b}

\newcommand{\Pgamma}[1]{P_\gga(#1)}
\newcommand{\normH}[1]{|#1|_H}

\newcommand{\Ldeu}{L^2(X,\gga)}

\newcommand{\sprod}[2]{\langle #1, #2 \rangle}

\newcommand\scal[2]{{\left\langle #1 ,#2\right\rangle}}

\def\diver{\mathop{\rm div}\nolimits}

\newcommand{\res}{\mathop{\hbox{\vrule height 7pt width .5pt depth 0pt
\vrule height .5pt width 6pt depth 0pt}}\nolimits}

\title[$BV$ functions in Wiener spaces]
{An introduction to $BV$ functions in Wiener spaces}

\author[Miranda-Novaga-Pallara]{M. Miranda jr, M. Novaga, D. Pallara}

\address{MM, Dipartimento di Matematica e Informatica, Universit\`a di Ferrara, via Machiavelli 35, 44121 Ferrara, Italy, 
{\rm michele.miranda@unife.it}
\\ 
MN, Dipartimento di Matematica, Universit\`a di Padova,\! via Trieste 63, 35121\! Padova, Italy,
{\rm novaga@math.unipd.it}
\\
DP, Dipartimento di Matematica e Fisica ``Ennio De Giorgi'', Universit\`a del Salento, P.O.B. 193, 
73100 Lecce, Italy, {\rm diego.pallara@unisalento.it}
}



\begin{document}

\begin{abstract}
We present the foundations of the theory of functions of bounded variation and sets
of finite perimeter in abstract Wiener spaces. 
\end{abstract}

\maketitle

\tableofcontents

\section{Introduction}\label{Intro}

This paper is an extended version of two talks given by the second and third author during
the summer school {\em Variational methods for evolving objects}. As both talks were concerned 
with some infinite dimensional analysis, we took the opportunity of this report  
to present the whole research area in a quite self-contained way, as it arises today. 
Indeed, even though geometric analysis on infinite dimensional spaces and the theory of $BV$ 
functions is presently an active research field and there are still many important open 
problems (some are presented in Section \ref{openproblems}), the foundations of the theory 
and some methods that have proved to be useful are 
rape enough as to be presented in an introductory paper. In particular, we think that our 
purpose fits into the general aim of a collection of lecture notes -- that of being useful
to students and young researchers who attended the summer school and could be interested 
in having an active part in further developments of the theory. 

Malliavin calculus is essentially a differential calculus in Wiener spaces and was initiated 
by P. Malliavin~\cite{mall76} in the seventies with the aim, among the others, of obtaining a probabilistic proof 
of H\"ormander hypoellipticity theorem. This quickly led to study connections to stochastic 
differential equations and applications in various fields in Mathematics and Physics, such as 
mathematical finance, statistical mechanics and hydrodynamics and the path approach to 
quantum theory or stationary phase estimation in stochastic oscillatory integrals with 
quadratic phase function. In general, solutions of SDEs are not continuous (and sometimes 
not even everywhere defined) functionals, hence the notion of weak derivative and Sobolev 
functional comes into play. Notice that there is no Sobolev embedding in the 
context of Malliavin calculus, which requires very little regularity. Looking at 
weak differentiation and the study of the behaviour of stochastic processes in domains 
leads immediately to the need for a good comprehension of integration by parts formulae, 
something that in the Euclidean case has been completely understood in the frameworks of 
geometric measure theory, sets with finite perimeter and more generally functions of bounded 
variation. This approach has been considered by Fukushima in \cite{fuk2000_1} and 
Fukushima-Hino in \cite{fuk2000_2}, where the first definition of $BV$ functions in infinite
dimensional spaces has been 
given, most likely inspired by a stochastic characterization of finite perimeter sets 
in finite dimension given by Fukushima in \cite{fuk99}, see Theorem 
\ref{fuk99result} below. In this paper we follow the integralgeometric approach to $BV$ 
functions developed in \cite{AMMP}, \cite{AMP}, \cite{AmbFig2}, \cite{AmbFigRuna}. 
Among the first applications of the theory, let us mention some results in a geometric 
vein in \cite{CasLunMirNov}, \cite{CasMirNov1} and in a probabilistic vein in \cite{trevisan}. 
On a more analytical perspective, some results are available on integral functionals, see  
\cite{CNV}, \cite{ChaGolNov}, and weak flows with Sobolev vector fields, see \cite{AmbFig1}. 
In this connection, the extension to $BV$ vector fields seems to require the analysis of fine 
properties of $BV$ functions and perimeters. \\
{\bf Acknowledgements.} The second and the third author are very grateful to the organisers 
of the summer school held in Sapporo in the summer of 2012, Professors L. Ambrosio, Y. Giga, 
P. Rybka and Y. Tonegawa, for the kind invitation and to the whole staff of the conference as well 
for their very pleasant stay in Japan. \\ This paper is partially supported 
by the Project ``Problemi di evoluzione e teoria geometrica della misura in spazi metrici''
of INdAM-GNAMPA. The second author acknowledges partial support by the Fondazione CaRiPaRo 
Project ``Nonlinear Partial differential Equations: models, analysis, and control-theoretic problems''.

\section{Preliminaries}\label{preliminaries}

As explained in the Introduction, motivations and possible 
applications of the theory we are going to present come from different areas, as well as the possible 
audience of the present notes. Indeed, it sits in the intersection between Calculus of Variations, 
Geometric Measure Theory, Functional Analysis, Stochastics and Mathematical Physics.  
Therefore, we have collected several prerequisites, divided in subsections, also with the purpose of 
fixing notation and basic results. Our aim is to introduce basic ideas and connections between the 
different perspectives, rather than giving precise and general results (this would take too much room). 
At the end of each subsection some general references for the sketched 
arguments are indicated. \\
When dealing with finite dimensional spaces $\R^d$, we always use Euclidean inner product 
$x\cdot y$ and norm $|x|$. Balls of radius $\varrho$ and centre $x$ in a Banach space are denoted 
by $B_\varrho(x)$, omitting the centre if $x=0$. The $\sigma$-algebra
of Borel sets in $X$ is denoted by $\B(X)$. Moreover,  we denote by 
$\|\cdot\|_X$ the norm in the Banach space $X$ and by $X^*$ the topological dual, 
with duality $\langle\cdot,\cdot\rangle$. 

\subsection{Measure theory}\label{Measuretheory}

In this subsection we briefly discuss a few properties of general measures with some details 
on {\em Gaussian measures} in finite and infinite dimensions. \\
A {\em measurable space} is a pair $(X,\F)$, where $X$ is a set and $\F$ a $\sigma$-algebra of 
subsets of $X$.  By {\em measure} on $(X,\F)$ we mean a countably additive function on $\F$ 
with values in a normed vector space; if a measure $\mu$ is given on $(X,\F)$, we say that 
$(X,\F,\mu)$ is a {\em measure space} (a {\em probability space} if $\mu$ is positive and 
$\mu(X)=1$) and omit $\F$ whenever it is clear from the context
or $\F=\B(X)$. For a measure $\mu$ with values in a normed vector space $V$ with norm $\|\cdot\|_V$ 
we define the total variation $|\mu|$ as the real valued positive measure
\begin{align*}
|\mu|(B)=\sup \Bigl\{
\sum_{j\in \N} \|\mu(B_j)\|_V : & B=\bigcup_{j\in \N} B_j,\ B_j\in\F,
\\
B_j\cap B_h=\emptyset \mbox{ for }j\neq h \Bigr\};
\end{align*}
the measure $\mu$ is said to be finite if $|\mu|(X)<+\infty$.
Given two measurable spaces $(X,\F)$ and $(Y,{\mathcal G})$, 
a measure $\mu$ on $X$ and a measurable function $f:X\to Y$ (i.e., such that $f^{-1}(B)\in\F$ 
for all $B\in{\mathcal G}$), the {\em push-forward} measure $\nu=f_\#\mu$ on $Y$ is defined by setting  
$\nu(B)=\mu(f^{-1}(B))$ for every $B\in{\mathcal G}$. Let us also recall that, given two measure 
spaces $(X_1,{\mathcal B}_1,\mu_1)$ and $(X_2,{\mathcal B}_2,\mu_2)$, the {\em product measure} 
$\mu_1\otimes\mu_2$ is defined on $X_1\times X_2$ by first defining the product $\sigma$-algebra
${\mathcal B}$ as that generated by $\{B_1\times B_2,\ B_1\in {\mathcal B}_1, B_2\in {\mathcal B}_2\}$ and 
then defining $\mu_1\otimes\mu_2$ as the unique measure on ${\mathcal B}$ such that 
$\mu_1\otimes\mu_2(B_1\times B_2)=\mu_1(B_1)\mu_2(B_2)$ for all pairs $B_j\in{\mathcal B}_j$. 
The construction generalises to the product of several spaces. 

In $\R^d$ we consider as reference measure either the Lebesgue measure $\LL^d$ or some absolutely 
continuous measure $\lambda=\rho\LL^d$ with nonnegative density $\rho$. The main 
examples among these are {\em Gaussian measures}. For $d=1$, these measures have densities $G$ given by 
\begin{equation}\label{defgaussian1dim}
G(x)=\frac{1}{\sqrt{2\pi q}} \exp\{-|x-a|^2/2q\}
\end{equation}
for some $a\in\R$ ({\em centre} or {\em mean}) and $q>0$ ({\em variance}). For $d>1$, a measure $\lambda$ on 
$\R^d$ is Gaussian if $f_\#\lambda$ is Gaussian on $\R$ for every linear function  
$f:\R^d\to\R$. Generalising \eqref{defgaussian1dim}, a Gaussian measure $\gamma$ on 
$\R^d$ is characterized by its centre $a=\int_{\R^d}x\, d\gamma$ and its {\em covariance matrix} 
$Q=(q_{hk})$ with 
\begin{equation}\label{qhk}
q_{hk} = \int_{\R^d} (x_h -a_h)\, (x_k-a_k)\ d\gamma(x),
\quad h,k\in\{1,\ldots, d\}  
\end{equation}
and is denoted $\NN(a,Q)$. A Gaussian measures $\gamma$ is {\em nondegenerate} if $\gamma=G\LL^d$ with 
$G(x)>0$ (equivalently, $Q$ positive definite) for all $x\in\R^d$, and is {\em standard} if 
\begin{equation}\label{standardgaussian}
G(x)= G_d(x)= (2\pi)^{-d/2}\exp\{-|x|^2/2\} ,
\end{equation}
i.e., $\gamma=\NN(0,{\rm Id})$. According to the preceding 
discussion on products, a standard Gaussian measure $\gamma_d$ on $\R^d=\R^k\times\R^m$ factors 
in the product of standard Gaussian measures $\gamma_d=\gamma_k\otimes\gamma_m$ for $k+m=d$. 
A measure $\gamma$ on a Banach space  
$(X,\B(X))$ is said Gaussian if $x^*_{\#}\gamma$ is Gaussian in $\R$ for every $x^*\in X^*$. 
In this case, the centre is defined as above by ({\em Bochner integral}, see \cite{boga})
\begin{equation}\label{center}
a=\int_{X}x\, d\gamma(x)
\end{equation}
and the covariance operator $Q\in {\mathscr L}(X^*,X)$ is a symmetric and positive operator uniquely 
determined by the relation, cf \eqref{qhk},
\begin{equation}\label{CovOp}
\scal{Qx^*}{y^*}=\int_X \scal{x-a}{x^*}\scal{x-a}{y^*}d\gamma(x),\ \forall x^*,y^*\in X^*.
\end{equation}
The fact that the operator $Q$ defined by \eqref{CovOp} is bounded is a consequence of Fernique's 
theorem (see e.g. \cite[Theorem 2.8.5]{boga}), asserting the existence of a positive $\beta>0$ such that
\begin{equation}\label{fernique}
\int_X \exp\{\beta\|x\|_X^2\}d\gamma(x)<\infty;
\end{equation}
indeed, $Q$ belongs to a special ideal of compact operators called {\em $\gamma$-Radonifying}. 
As above, we write $\gamma=\NN(a,Q)$ and we say that $\gamma$ is nondegenerate if ${\rm Ker}\, Q=\{0\}$. 
Notice that the Dirac measure at $x_0$ is considered as the (fully 
degenerate) Gaussian measure with centre $x_0$ and covariance $Q=0$. 
For the arguments of the present subsection we refer to \cite{boga}, \cite {bogamt}. 

\subsection{Geometric measure theory}\label{gmt}

A general class of (non absolutely continuous) measures of interest in the sequel is that of 
{\em Hausdorff measures}, which we briefly discuss here, together with the related notions
of {\em rectifiable set} and {\em approximate tangent space}. \\
The measure ${\mathcal H}^s$, $0<s<\infty$, is defined in a general metric space by 
\begin{equation}\label{defHausdorffmeas}
{\mathcal H}^s(B) = \frac{\omega_s}{2^s}\sup_{\delta>0} \inf\Bigl\{\sum_{j=1}^\infty ({\rm diam}\, B_j)^s,\ 
B\subset \bigcup_{j=1}^\infty B_j,\ {\rm diam}\, B_j < \delta\Bigr\},
\end{equation}
where, using Euler's $\Gamma$ function, $\omega_s= \Gamma(1/2)^s/\Gamma(s/2+1)$ ($=\LL^d(B_1)$ 
if $s=d\in\N$) is a normalising constant and the infimum runs along all the countable coverings. 
Beside the Hausdorff measures, it is useful to introduce the {\em Minkowski content}, which provides 
a more elementary, though less efficient, way of measuring ``thin'' sets. Given a closed set 
$C\subset\R^d$ and an integer $s$ between $0$ and $d$, the idea is to look 
at the rate of convergence to $0$ of $\varrho\mapsto \LL^d\left(I_\varrho(C)\right)$ as
$\varrho\downarrow 0$, where $I_\varrho(C)$ denotes the open $\varrho$-neighbourhood of $C$.  
In general, given a closed set $C\subset\R^d$, the upper and lower $s$-dimensional Minkowski 
contents ${\mathcal M}^{*s}(C)$, ${\mathcal M}_*^s(C)$ are defined by
\begin{equation}\label{defMinkowski}
\begin{array}{l}
\displaystyle{
{\mathcal M}^{*s}(C)=\limsup_{\varrho\downarrow 0}
\frac{\LL^d(I_\varrho(C))}{\omega_{d-s}\varrho^{N-s}},}
\\ \ \\ 
\displaystyle{
{\mathcal M}^s_*(C)=\liminf_{\varrho\downarrow 0}
\frac{\LL^d(I_\varrho(C))}{\omega_{d-s}\varrho^{N-s}},}
\end{array}
\end{equation}
respectively. If ${\mathcal M}^{*s}(S)={\mathcal M}_*^s(C)$, their common value is denoted by 
${\mathcal M}^s(C)$ (Minkowski content of $C$) and we say that $C$ admits Minkowski content. 
Unlike the Hausdorff measures, the Minkowski content 
is not subadditive. Nevertheless, in some important cases the two procedures give the same result. 
We compare later the Hausdorff measures and the Minkowski contents. 

The natural regularity category in geometric measure theory is that of {\em Lipschitz continuous 
functions}. Let us recall (Rademacher theorem) that a Lipschitz function defined on $\R^d$ with 
values in a finite dimensional vector space is differentiable $\LL^d$-a.e. (the differentiability 
properties of Lipschitz functions defined on infinite dimensional vector spaces is a much more delicate 
issue, see \cite{boga}, \cite{preiss}). For $s$ integer between 0 and $d$, we say that a ${\mathcal H}^s$ 
measurable set $B\subset \R^d$ is {\em countably $s$-rectifiable} if there are countably many Lipschitz
functions $f_j:\R^s\to\R^d$ such that
\begin{equation}\label{rectset}
B\subset\bigcup_{j=1}^\infty f_j(\R^s).
\end{equation}
We say that $B$ is {\em countably ${\mathcal H}^s$-rectifiable} if there
are countably many Lipschitz functions $f_j:\R^s\to\R^d$ such that
\begin{equation}\label{countrectset}
{\mathcal H}^s\Bigl(B\setminus\bigcup_{j=0}^\infty f_j(\R^s)\Bigr)=0.
\end{equation}
Finally, we say that $B$ is {\em ${\mathcal H}^s$-rectifiable} if $B$ is countably 
${\mathcal H}^s$-rectifiable and ${\mathcal H}^s(B)<\infty$. All these classes of sets are stable 
under Lipschitz mapping. Notice that countable ${\mathcal H}^s$-rectifiability is equivalent to 
the seemingly stronger requirement that ${\mathcal H}^s$-almost all of the set can be covered by 
a sequence of Lipschitz $s$-graphs. Notice that if the admissible coverings in \eqref{defHausdorffmeas}
are made only by balls we get the {\em spherical} Hausdorff measure ${\mathcal S}^s$. The measures 
${\mathcal H}^s$ and ${\mathcal S}^s$ are comparable in the sense that
\[
{\mathcal H}^s \leq {\mathcal S}^s \leq 2^s {\mathcal S}^s
\]
and coincide on ${\mathcal H}^s$-rectifiable sets. However,  
an important difference between ${\mathcal H}^s$ and ${\mathcal S}^s$ measures is relevant in 
Subsection \ref{SubsHausdorff}, where Hausdorff measures are discussed in the infinite dimensional 
setting, see Lemma \ref{lmono}. Analogously, the Hausdorff measure coincide with the Minkowski  
content on rectifiable sets. Even though rectifiable sets can be very irregular from the 
point of view of classical analysis, nevertheless they enjoy useful properties from the point of view of 
geometric measure theory. Indeed, for ${\mathcal H}^s$-a.e point $x$ of a countably 
${\mathcal H}^s$-rectifiable set $B$ there exists an $s$-dimensional subspace  
$S$ ({\em approximate tangent space}) such that 
\begin{equation}\label{apptgspace}
\lim_{\varrho\to 0} \int_{\frac{B-x}{\varrho}} \varphi\, d{\mathcal H}^s = 
\int_S \varphi\, d{\mathcal H}^s \qquad \forall\ \varphi \in C_c(\R^d). 
\end{equation}
If $s=d-1$ an {\em approximate unit normal vector} $\nu(x)$ to $B$ at $x$ is defined (up to the 
sign) as the unit vector normal to $S$. In the same vein, we say that a function 
$u\in L^1_{\rm loc}(\R^d,\R^k)$ admits an {\em approximate limit} at $x_0$, if
there is $z\in\R^k$ such that 
\begin{equation}\label{defapplim}
\lim_{\varrho\to 0} \frac{1}{\omega_d\varrho^d} \int_{B_\varrho(x_0)} |u(x)-z|\, dx = 0 
\end{equation}
($z=\aplim_{x\to x_0}u(x)$ for short) and in this case we say that $u$ is {\em approximately 
continuous at $x_0$} if $x_0$ is a Lebesgue point of $u$ and \eqref{defapplim} holds with $z=u(x_0)$. 
Analogously, if $u$ is approximately continuous at $x_0$ we say that $u$ is {\em approximately 
differentiable at $x_0$} if there is a linear map $L:\R^d\to\R^k$ such that
\begin{equation}\label{defappdiff}
\aplim_{x\to x_0}\frac{u(x)-u(x_0)-L(x-x_0)}{|x-x_0|}=0.
\end{equation}
For the arguments of the present subsection we refer to \cite{afp}, \cite {federer}. 

\subsection{Stochastic analysis}\label{stochastics}

Let a probability space $(\Omega,\F,\PP)$ be given. If $(X,{\mathcal B})$ is a measurable 
space, a measurable function $\xi:\Omega\to X$ is called an {\em $X$-valued random variable} 
(r.v. for short) and its {\em law} is the push-forward measure of $\PP$ under $\xi$, i.e., 
$\xi_{\#}\PP(B)=\PP(\xi^{-1}(B))$, $B\in \B$. If $\xi\in L^1(\Omega,\PP)$ we define its 
{\em expectation} by $\E[\xi]=\int_\Omega\xi\, d\PP$; if $\xi\in L^2(\Omega,\PP)$ we define 
its {\em variance} by $Var(\xi) = \E[\xi - \E[\xi]] = \E[\xi^2]-\E^2[\xi]$ and for 
$\xi,\eta\in L^2(\Omega,\PP)$ we define the {\em covariance} by 
${\rm cov}(\xi,\eta) = \E[\xi\eta]-\E[\xi]\E[\eta]$. Given a sub $\sigma$-algebra 
${\mathcal G}\subset\F$, the {\em conditional expectation of a summable $\xi$ given ${\mathcal G}$} 
is the unique ${\mathcal G}$-measurable random variable $\eta=\E(\xi|{\mathcal G})$ such that 
$\int_B\xi\, d\PP= \int_B\eta\, d\PP$ for all $B\in {\mathcal G}$. Given $N$ random variables 
$\xi_j:\Omega\to X_j$, they are {\em independent} if for every $A_j\subset X_j$, setting 
$B_j=\{\omega\in\Omega:\ \xi_j(\omega)\in A_j\}$, $\PP(B_1\cap\cdots\cap B_N)=\PP(B_1)\cdots \PP(B_N)$, 
or, equivalently, if the law of the r.v. $\xi:\Omega\to X=X_1\times\cdots\times X_N$ whose 
components are the $\xi_j$ is the product measure of the laws of the $\xi_j$ on $X$. A random variable 
is Gaussian if its law is a Gaussian measure. 

An $X$-valued {\em continuous stochastic process $\xi$ on $[0,\infty)$} is the assignment, for 
$t\in[0,\infty)$, of a family of random variables $\xi_t:(\Omega,\F,\PP)\to X$. An
increasing family of sub $\sigma$-algebras $\F_t\subset \F$ is called a {\em filtration};
a process $\xi$ is said {\em adapted} to a given filtration $\F_t$ if $\xi_t$
is $\F_t$-measurable for every $t$.
If the filtration is not explicitly assigned, the {\em natural filtration} is understood, i.e., $\F_t$ is 
the smallest $\sigma$-algebra such that $\xi_s$ is measurable for all $s\leq t$, $s\in I$. 
If $\xi_t$ is an adapted process, summable for every $t$ and $\E(\xi_t|\F_s)=\xi_s$ 
for all $s\leq t$, the process $\xi$ is a {\em martingale}. 
Due to the dependence of $\xi_t(\omega)$ on two variables, we may think of $\omega\mapsto\xi_t(\omega)$, 
for fixed $t$, as a family of r.v. defined on $\Omega$, or as $t\mapsto\xi_t(\omega)$, for $\omega$ 
fixed, as a set of {\em trajectories}. A real stochastic process on an interval $I$ defines the  
{\em distribution functions}
\[
F_{t_1\cdots t_n}(x_1,\ldots, x_n)=\PP[\xi_{t_1}<x_1,\ldots,\xi_{t_n}<x_n],
\ 0\leq t_1<\ldots<t_n<\infty, 
\]
called {\em finite-dimensional joint distributions}. In general, $F(x)$ is said to 
be a distribution function if it is increasing with respect to all the $x_k$ 
variables, left-continuous, $F(x_1,\ldots, x_n)\to 0$ if some $x_k\to-\infty$, 
$F(x_1,\ldots, x_n)\to 1$ if all $x_k\to+\infty$ and for any intervals 
$I_k=[a_k,b_k)$, $1\leq k\leq n$, the inequality
\[
\Delta_{I_1}\cdots \Delta_{I_n} F(x_1,\ldots,x_n) \geq 0
\]
holds, where $\Delta_{I_k}F(x)=F(x_1,\ldots,b_k,\ldots,x_n)-F(x_1,\ldots,a_k,\ldots,x_n)$. 
A (remarkable) result of Kolmogorov's states that, given a sequence $F_n(x_1,\ldots,x_n)$ of 
distribution functions, there is always a stochastic process whose distribution functions are 
the given ones, provided the (necessary) {\em consistency condition} 
\[
\lim_{x_n\to+\infty}F_n(x) = F_{n-1}(x_1,\ldots,x_{n-1})
\]
holds. A stochastic process on $[0,\infty)$ is {\em stationary} if its distribution function is 
invariant under translations on time, i.e.,
\[
F_{t_1+h\cdots t_n+h}(x_1,\ldots, x_n) = F_{t_1\cdots t_n}(x_1,\ldots, x_n)\qquad \forall\ h\geq 0. 
\]
Given a filtration $\F_t,\ t\in I$, a random variable $\tau:\Omega\to I=[0,+\infty]$ 
is a {\em stopping time} if $\{\tau\leq t\}\in \F_t$ for all $t\in I$. Accordingly, a process 
$\xi$ is a {\em local martingale} if there is an increasing sequence of stopping times 
$\tau_n\to +\infty$ such that $(\xi_{t\wedge\tau_n})$ is a martingale for every $n\in\N$. 

A particular class of processes which is relevant for our purposes is that of {\em Markov processes}. 
Let us start from the notion of {\em time-homogeneous Markov transition function}, i.e., a function 
$p(t,x,B)$, $t\in [0,\infty)$, $x\in X$, $B\in {\mathcal B}$, which is measurable with respect to $x$, 
is a probability measure on $(X,{\mathcal B})$ with respect to $B$ (we also write $p(t,x,dy)$ to stress
the last property) and verifies the {\em Chapman-Kolmogorov} equation
\begin{equation}\label{chapmankol}
p(t,x,B) = \int_X p(t-s,y,B)p(s,x,dy) ,\qquad \forall\ 0\leq s \leq t.
\end{equation}
Given a transition function $p$ as above and a probability distribution $\mu$ on $(X,{\mathcal B})$, 
there is a stochastic process $\xi$ such that the law of $\xi_0$ is $\mu$ and
$\PP(\xi_t\in B|\F_s)=p(t,\xi_s,B)$ for all $0\leq s\leq t$ and it is called 
{\em Markov process associated with $p$ with initial law $\mu$}. The initial law $\mu$ is
{\em invariant} with respect to the process (see also next Subsection) if 
\begin{equation}\label{invariant1}
\mu(B) = \int_X p(t,y,B)\, \mu(dy) ,\qquad \forall\ t\geq 0,\ B\in {\mathcal B}.
\end{equation}
An $\R^d$ valued $Q$-{\em Brownian motion} starting from $a$ or {\em Wiener process} 
$B_t$ is a stochastic process such that $B_0=a\in\R^d$, for $\PP$-a.e. $\omega\in\Omega$ the trajectories 
$t\mapsto B_t(\omega)$ are continuous, for every $0\leq s<t$ the difference $B_t-B_s$ is a 
Gaussian random variable with centre $a$ and covariance $(t-s)Q$, i.e., $\NN(a,(t-s)Q)$ and for every 
$0\leq t_1<\ldots<t_n$ 
the random variables $B_{t_2}-B_{t_1},\ldots,B_{t_n}-B_{t_{n-1}}$ are independent. This in partucular implies that
the Brownian motion is a martingale since the independence of $B_t-B_s$ from $B_s$ implies that
$B_t-B_s$ is independent from ${\mathcal F}_s$, that is
\[
\E(B_t-B_s|\F_s)=0.
\]
According to the quoted Kolmogorov theorem, Brownian motions exist. Notice that a Brownian motion is 
a Markov process whose transition function is Gaussian, $p(t,x,dx)=\NN(x,tQ)$. Moreover, as we have 
already observed, any Brownian motion has a continuous version and is a martingale; in the sequel we 
always assume that the continuous version has been selected. A Brownian motion is {\em standard} (or 
{\em normalised}) if $a=0$ and $Q={\rm Id}$. 

The {\em It\^o integral} with respect to a given real Brownian motion $B_t$, whose (completed) natural
filtration we denote by $\F_t$, can be defined through suitable Riemannian sums, even though 
the usual Stiltjes approach cannot be pursued, due to the fact that $B_t$ has not bounded variation 
in time. Let $\xi_t,\ t\in [0,T]$ be an adapted continuous {\em simple} process, i.e.,  such 
that there are a partition $0=t_0<t_1< \cdots < t_N\leq T$ and $\F_{t_{j-1}}$-measurable 
r.v. $\xi_j,\ j=1,\ldots, N$,  for which
\[
\xi_t(\omega) = \sum_{j=1}^N \xi_j(\omega)\chi_{[t_{j-1},t_j)}(t).
\]
For such a process, define 
\[ 
\int_0^T\xi_t\, dB_t = \sum_{j=1}^N \xi_j(B_{t_j}-B_{t_{j-1}}). 
\] 
As a consequence of the independence of the increments of the Brownian motion, 
we get the {\em It\^o isometry}
\begin{equation}\label{ItoIsom}
\PP\Bigl(\int_0^T\xi_t\, dB_t \cdot \int_0^T\eta_t\, dB_t\Bigr) = \int_0^T \xi_t \eta_t\, dt
\end{equation}
for every $\xi,\eta$ as above. The It\^o isometry extends to $\R^d$ valued processes and 
Brownian motion in an obvious way.
Thanks to the It\^o isometry and the fact that every adapted process $\xi$ such that $\PP(\int_0^T|\xi_s|^2ds<\infty)=1$ 
can be approximated by elementary processes, it is possible to extend the stochastic integral 
to the described class of processes, or to processes defined for $0\leq t<\infty$ such that 
the finiteness condition holds for every $T>0$. Notice that the stochastic integral is, in turn, 
a random variable. It can be proved as well that the function $t\mapsto\int_0^t \xi_s\, dB_s$ is 
continuous $\PP$-a.s.  

The stochastic integral allows for a rigorous theory of {\em stochastic differential equations}, SDEs for 
short, which are intuitively dynamical systems perturbed by noise. We deal here only with {\em autonomous 
SDEs} on $\R^d$, assuming that the noise is given in terms of a Brownian motion. Something more in the 
Wiener space will be added in Subsection \ref{OUsection} in connection with the OrnsteinUhlenbeck process. 
In the present case the Cauchy problem can be written (at least formally) as 
\begin{equation}\label{SDE}
d\xi_t=A(\xi_t)dt+\sigma(\xi_t)dB_t,\qquad \xi_0\ \ {\rm given\ r.v.}, 
\end{equation}
where $B_t$ is a Brownian motion, $\sigma$ and $A$ are the {\em diffusion} and {\em drift}
term, respectively. The meaning of \eqref{SDE} is that the process $\xi$ is a solution if 
\[
\xi_t = \xi_0 + \int_0^t A(\xi_s)\, ds + \int_0^t \sigma(\xi_s)\, dB_s.
\]
Under general hypotheses a unique solution exists and is a continuous Markov process. 
Presenting a general theory goes far from the aim of this short presentation; detailed 
results are discussed on concrete cases. For the arguments of the present subsection we 
refer to \cite{boga}, \cite {friedman}, \cite{malliavin}. 

\subsection{Semigroup theory}\label{semigroups}

The theory of one-parameter semigroups of linear operators in Banach spaces was born 
as a general method to solve autonomous evolution equations, has been widely studied 
and is very rich of abstract results and applications. We need very few basic results,  
and the main point which is worth discussing here is the link between semigroups as a 
tool for solving linear parabolic partial differential equations and the related 
stochastic differential equations, as explained at the end of this subsection. First, 
we say that $(S_t)_{t\geq 0}$ is a {\em semigroup of linear operators} on a Banach 
space $E$ if $S_t\in {\mathscr L}(E)$, i.e., $S_t$ is a bounded linear operator on 
$E$ for every $t\geq 0$, $S_0={\rm Id}$, $S_{t+s}=S_t\circ S_s$; if $t\mapsto S_t f$ 
is norm continuous for every $x\in E$ then $S_t$ is said to be $C_0$ (or strongly 
continuous). If $S_t$ is strongly continuous then, setting 
$\omega_0=\inf_{t\geq 0}\frac{1}{t}\log\|S_t\|_{{\mathscr L}(E)}$, for every $\varepsilon$ there is 
$M_\varepsilon\geq 1$ such that $\|S_t\|_{{\mathscr L}(E)} \leq M_\varepsilon e^{(\omega_0+\varepsilon) t}$ 
for all $t\geq 0$.  A semigroup defined on $E=C_b(X)$ (the space of bounded continuous functions on a 
Banach space $X$) is {\em Feller} if $S_tf\in C_b(X)$ for all $f\in C_b(X)$ and is {\em strong Feller} 
if $S_tf\in C_b(X)$ for all $f\in B_b(X)$ (the space of bounded Borel functions). 
A {\em Markov semigroup} is a semigroup $S_t$ on $C_b(X)$ such that $S_t\one = \one$, 
$\|S_t\|_{{\mathscr L}(E)}\leq 1$ for every $t\geq 0$, and $S_tf \geq 0$ for every $f\geq 0$ and $t>0$ 
(here $\one$ is the constant function with value 1). Given a time homogeneous Markov transition function 
$p$ and the associated process $\xi_t^x$ starting at $x$ (which means that the law of $\xi_0$ is $\delta_x$), 
the family of operators
\begin{equation}\label{Markovsgrp}
S_tf(x) = \int_{X} f(y) p(t,x,dy) = \E[f(\xi_t^x)] ,\qquad x\in X,
\end{equation}
due to \eqref{chapmankol}, is a {\em Markov semigroup}. Notice that $S_t$ can be extended to $B_b(X)$. 
With each semigroup
it is possible to associate a {\em generator}, i.e., a linear closed operator $(L,D(L))$ such that
$Lf=\lim_{t\to 0}(S_tf-f)/t$, $f$ in the {\em domain} $D(L)\subset E$. Here the limit is in the norm sense 
if $S_t$ is strongly continuous or can be in weaker senses (uniform convergence on bounded or compact sets 
or even pointwise with bounds on the sup norm) in the case of Markov semigroups. 
We are mainly interested in the case where $p$ comes from a process which solves a SDE \eqref{SDE} on a 
Banach space $X$. In this case, $L$ is a linear elliptic operator given by 
$L=-\frac{1}{2}{\rm Tr}[\sigma\sigma^*D^2]+\scal{Ax}{\nabla}$, at least on suitable smooth functions,  
giving rise to the {\em Kolmogorov backward parabolic operator} $\partial_t-L$. Under suitable conditions, 
the solution of the Cauchy problem $\partial_tu-Lu=0$, $u(0)=f\in C_b(X)$ will be given by $u(t)=S_tf$. 
In this setting, the trajectories of the Markov process play a role analogous to that of the 
characteristic curves in a hyperbolic problem. Finally, we introduce the notion of {\em invariant
measure} associated with the semigroup $S_t$, i.e., a probability measure $\mu$ on $X$ such that 
\[
\int_X S_tf (x)\, d\mu(x) = \int_X f(x)\, d\mu(x), \qquad f\in C_b(X).
\]
The meaning of the above equality is that the distribution $\mu$ is invariant under the 
flow described by equation \eqref{SDE}, see \eqref{invariant1}. Typically, if $\mu_t$ is the 
law of $\xi_t$ and the weak limit $\mu=\lim_{t\to\infty}\mu_t$ exists, then $\mu$ is invariant
and the semigroup $S_t$ extends to a $C_0$ semigroup in all the $L^p(X,\mu)$ spaces, $1\leq p < \infty$. 
For the arguments of the present subsection we refer to \cite{boga}, \cite {friedman}. 

\subsection{Dirichlet forms}\label{dirichletforms}

In this subsection we collect a few notions on Dirichlet forms, confining to what we need 
in Theorem \ref{fuk99result}, and to show some further 
connections between the various areas we are quickly touching. 

Given a $\sigma$-finite measure space $(X,\mu)$ consider the Hilbert space $L^2(X,\mu)$ with the 
inner product $[u,v]$. A functional $\EE: D(\EE)\times D(\EE)\to\R$ is a {\em Dirichlet form} if it is 
\begin{enumerate}
\item {\em bilinear}: $\EE(u+v,w)=\EE(u,w)+\EE(v,w)$,\ $\EE(\alpha u,v)=\alpha\EE(u,v)$
for all $u,v,w\in L^2(X,\mu)$, $\alpha\in \R$;
\item {\em nonnegative}: $\EE(u,u)\geq 0$ for all $u\in L^2(X,\mu)$;
\item {\em closed}: $D(\EE)$ is complete with respect to the metric induced by the inner
product $\EE(u,v)+[u,v]$,\ $u,v\in D(\EE)$;
\item {\em Markovian}: if $u\in D(\EE)$ then $v:=(0\vee u)\wedge 1\in D(\EE)$ and 
$\EE(v,v)\leq \EE(u,u)$.
\end{enumerate}
A Dirichlet form $\EE$ is {\em symmetric} if $\EE(u,v)=\EE(v,u)$ for all $u,v\in L^2(X,\mu)$ and 
is {\em local} if $\EE(u,v)=0$ whenever $u,v\in D(\EE)$ have disjoint compact supports. 
The subspace $D(\EE)$ of $L^2(X,\mu)$ is called the {\em domain} of the form $\EE$.
 
Dirichlet forms are strictly connected with Markov semigroups and processes. First, 
notice that a nonnegative operator $L$ can be associated with any Dirichlet form 
as shown in the following theorem of Kato's. 

\begin{theorem}\label{tKato}
There is a one-to-one correspondence between closed symmetric forms and 
nonnegative self-adjoint operators given by
\[
u\in D(L)\ \Leftrightarrow \ \exists\, f\in L^2(X,\mu):\ \ 
\EE(u,v)=[f,v]\ \forall\, v\in D(\EE),\ Lu:=f. 
\]
Moreover, $D(\EE)=D(\sqrt{L})$ and the operator $(-L,D(L))$ is the generator of a strongly 
continuous Markov semigroup $S_t$ of self-adjoint operators. 
\end{theorem}

According to the discussion in the preceding Subsection and the above Theorem, it is possible 
to associate with a Markov process, beside a Markov semigroup, a Dirichlet form. Of course, 
not all the Markov processes give raise to a Dirichlet form. Moreover, the transition function 
must be {\em symmetric}, i.e., such that $p(x,y,B)=p(y,x,B)$ for all $x,y\in X$ and $B\in \B(X)$
in order to get a symmetric Dirichlet form and if the process has continuous trajectories then 
the associated form is local. 

Viceversa, given a regular Dirichlet form, there is a unique (in a suitable sense) Markov process  
whose Dirichlet form is the given one. Let us now discuss two key examples that will play a relevant 
role in the sequel. 

\begin{example}\label{Neumann}
{\rm Let $D\subset\R^d$ be open and bounded with Lipschitz continuous boundary, 
and define the Dirichlet form on $L^2(D)$ by 
\[
\EE(u,v)=\int_{D} \nabla u\cdot\nabla v\ dx,
\]
for $u,v\in D(\EE)=W^{1,2}(D)$. The operator $L$ defined as in Theorem 
\ref{tKato} is the Neumann Laplacean, i.e., 
\[
L=-\Delta,\qquad D(L)=\{u\in H^{2,2}(D):\ \partial_\nu u=0\ {\rm on\ }\partial D\}, 
\]
where $\partial_\nu$ denotes the differentiation with respect the normal direction. 
Then, $(-L,D(L))$ is the generator of a strongly 
continuous Markov semigroup on $L^2(D)$ and the related Markov process is the {\em reflecting 
Brownian motion in} $D$. 
}\end{example}

\begin{example}\label{OUexample}
{\rm Let $\gamma=G_d\LL^d$ be the standard Gaussian measure. Define the Dirichlet form 
$\EE$ on $L^2(\R^d,\gamma)$ by 
\[
\EE(u,v)=\int_{\R^d} \nabla u\cdot\nabla v\ d\gamma,
\]
$u,v\in D(\EE)=W^{1,2}(\R^d,\gamma)=\{u\in W^{1,2}_{\rm loc}(\R^d)\!: u, |\nabla u|\in L^2(\R^d,\gamma)\}$.   
The operator $L$ defined as in Theorem \ref{tKato} is the Ornstein-Uhlenbeck operator  defined 
on smooth functions by $L=-\Delta+x\cdot\nabla$ and $D(L)= W^{2,2}(\R^d,\gamma)$, $(-L,D(L))$ 
is the generator of the strongly continuous Markov semigroup $T_t$ on $L^2(\R^d,\gamma)$ defined in 
\eqref{OUsgrpinRd} and the related Markov process is the {\em Ornstein-Uhlenbeck process in} $\R^d$
given by \eqref{OUprocess} below. Moreover, $\gamma$ is the invariant measure of $T_t$. 
}\end{example}
For the arguments of the present subsection we refer to \cite{Fot94}, \cite{MaRockner}. 

\section{$BV$ functions in the finite-dimensional case}\label{secfinitedim}

In this section we present the main properties of $BV$ functions in $\R^d$. 
In order to pave the way to the generalisations to Wiener spaces, we discuss now
at the same time the case when the reference measure is the 
Lebesgue one or the finite dimensional standard Gaussian measure. Of course, $BV$ functions 
with general densities can be studied, but this is not of our concern here. Standard Gaussian 
measures have regular and non-degenerate densities, hence there is no basic difference at the 
level of {\em local} properties of $BV$ functions, which are basically the 
same in the two cases. Instead, the {\em global} properties are different, 
due to the very different behaviour of the densities at infinity. Let us start 
from the classical case. There are various ways of defining $BV$ 
functions on $\R^d$, which are useful in different contexts. 

\begin{theorem}\label{findimeq}
Let $u\in L^1(\R^d)$. The following are equivalent:
\begin{itemize}
\item[{\bf 1}] there exist real finite measures $\mu_j,\ j=1,\ldots, d,$ on $\R^d$ such that 
\begin{equation}\label{intbyparts}
\int_{\R^d} uD_j\phi dx= - \int_{\R^d} \phi d\mu_j,\ \ \forall\phi\in C_c^1(\R^d), 
\end{equation}
i.e., the distributional gradient $Du=\mu$ is an $\R^d$-valued measure with finite total 
variation $|Du|(\R^d)$;
\item[{\bf 2}] the quantity 
\begin{equation*}
V(u)=\sup\Big\{\int_{\R^d} u\diver \phi dx:\phi\in C_c^1(\R^d,\R^d),
\|\phi\|_{\infty}\leq 1\Big\}
\end{equation*}
is finite;
\item[{\bf 3}] the quantity
\begin{equation*}
L(u)=\inf\Big\{\liminf_{h\to \infty}\int_{\R^d}|\nabla u_h|\, dx:
u_h\in {\rm Lip}(\R^d),\ u_h\stackrel{L^1}{\rightarrow}u\Big\}
\end{equation*}
is finite;
\item[{\bf 4}] if $(W_t)_{t\geq 0}$ denotes the heat semigroup in $\R^d$, then 
\begin{equation*}
{\mathscr W}[u]=\lim_{t\to 0} \int_{\R^d} |\nabla W_tu|\, dx <\infty.
\end{equation*}
\end{itemize}
Moreover, $|Du|(\R^d)=V(u)=L(u)={\mathscr W}[u]$.
\end{theorem}
If one of (hence all) the conditions in Theorem \ref{findimeq} holds, we 
say that $u\in BV(\R^d)$. The statement above is well known, a sketch of
its proof, with more references, can be found in \cite{AMMPPhysicad}. 
We observe that in {\bf 3} we may replace Lipschitz functions with functions in $W^{1,1}(\R^d)$.
The translation of the above result in the case of a standard Gaussian 
measure $\gamma=\NN(0,{\rm Id})=G_d\LL^d$ is an easy matter, taking into account that the 
integration by parts formula has to be modified because the density of $\gamma$ is 
not constant and reads
\begin{equation}\label{gaussintbyparts}
\int_{\R^d} u(x) D_jv(x) \, d\gamma(x) = 
-\int_{\R^d} \bigl[ v(x) D_ju(x) - x_j u(x) v(x)\bigr] \, d\gamma(x) .
\end{equation}
Hence, $BV(\R^d,\gamma)$ functions 
and the weighted total variation measure $|D_\gamma u|$ can be defined, for 
$u\in L^1(\R^d,\gamma)$, as in the above Theorem, according to the following suggestions:
\begin{enumerate}
\item replace the measure $dx$ with $d\gamma$ everywhere;
\item in {\bf 1}, replace $D_j\phi(x)$ with 
$D_j^*\phi(x)=D_j\phi(x)-x_j\phi(x)$;
\item in {\bf 2}, replace ${\rm div}\, \phi$ with $\sum_{j=1}^dD_j^*\phi_j$;
\item in {\bf 4}, replace the heat semigroup $W_t$ with the {\em Ornstein-Uhlen\-beck semigroup}
\begin{equation}\label{OUsgrpinRd}
\begin{array}{ll}
T_tu(x)&=\displaystyle{\int_{\R^d} u(e^{-t}x+\sqrt{1-e^{-2t}}y)\, d\gamma(y) }
\\
&\displaystyle{=(2\pi)^{-d/2}\int_{\R^d} u(e^{-t}x+\sqrt{1-e^{-2t}}y)e^{-|y|^2/2}\, dy}
\\
&\displaystyle{=(2\pi(1-e^{-2t}))^{-d/2} \int_{\R^d} u(y) e^{-|y-e^{-t}x|^2/2\sqrt{1-e^{-2t}}}dy}
\end{array}
\end{equation}
which plays a fundamental role in the infinite-dimensional analysis.  
\end{enumerate}
Using Dirichlet forms, a further characterization of $BV$ functions can be given in the Gaussian setting. 
Indeed, given $u\in L^1(\R^d,\gamma)$, for $j=1,\ldots, d$ the linear projections $x_j^*$ belong to the 
domain of the form 
\[
\EE^u(w,v)=\int_{\R^d} \nabla w\cdot\nabla v\ u\ d\gamma
\] 
and $u\in BV(\R^d,\gamma)$ if and only if there is $C>0$ such that 
\begin{equation}\label{BVconDirform}
|\EE^u(x^*_j,v)| \leq C \|v\|_\infty\qquad \forall\ v\in C^1_b(\R^d). 
\end{equation}
Notice that both the heat and the Ornstein-Uhlenbeck semigroups are Markov semigroups whose 
transition functions in $\R^d$ in the sense of \eqref{Markovsgrp} are given by
\begin{align} \nonumber 
p(t,x,dy) &= G(t,x,y)dy, 
\\ \label{heatkernel}
G(t,x,y) & = \frac{1}{t^{d/2}} G_d\Bigl(\frac{x-y}{\sqrt{t}}\Bigr) = 
\frac{1}{(2\pi t)^{-d/2}}\exp\Bigl\{\frac{|x-y|^2}{2t}\Bigr\}
\\ \nonumber 
p(t,x,dy) &= \psi(t,x,y)dy,
\\ \label{OUkernel} 
\psi(t,x,y) &= 
(2\pi(1-e^{-2t}))^{-d/2} \exp\Bigl\{-\frac{|y-e^{-t}x|^2}{2\sqrt{1-e^{-2t}}}\Bigr\}. 
\end{align}
The only non trivial point is (4), which is discussed in detail in the Wiener case. For the 
moment, as discussed also in Subsection \ref{dirichletforms} and in particular in Example \ref{OUexample}, 
let us only point out that the infinitesimal generator of $T_t$ is the operator defined on smooth 
functions by the expression 
\[
-Lu(x)=\Delta u(x)- x \cdot \nabla u(x) 
\]
and that $\gamma$ turns out to be the invariant measure associated with $T_t$. The
semigroup $T_t$ is related to the Ornstein-Uhlenbeck process 
\begin{equation}\label{OUprocess}
\xi_t= e^{-t/2}\xi_0 + \int_0^t e^{(s-t)/2}\, dB_s, 
\end{equation}
solution of the {\em Langevin SDE} 
\begin{equation}\label{Langevin}
d\xi_t=-\frac{1}{2} \xi_t\, dt + dB_t.  
\end{equation}
From this point of view, let us recall that the generator of $W_t$ is the Laplace operator, 
and that the Lebesgue measure is invariant under the heat flow (this does not fit completely
into the theory of invariant measures, as $\LL^d$ is not finite). 

Differently from the Sobolev case, $BV$ functions are allowed to be discontinuous along hypersurfaces, 
and indeed {\em characteristic functions} $\chi_E$ may belong to $BV$. If $E\subset\R^d$ and 
$|D\chi_E|(\R^d)$ is finite, we say that $E$ is a set with finite perimeter, and use the notation 
$P(E)$ (perimeter of $E$) for the total variation of the measure $D\chi_E$ and write $P(E,\cdot)$ 
for $|D\chi_E|(\cdot)$. Analogously, we set $P_\gamma(E)$ and $P_\gamma(E,\cdot)$ in the Gaussian
case. The study of structure of sets with finite perimeter is important on its own, but also because 
it gives information on general $BV$ functions, through the {\em coarea formula}: if $u\in BV(\R^d)$, 
then $P(\{u>t\})$ is finite for a.e. $t\in\R$ and for every $B\in \B(\R^d)$ the following equality holds:
\begin{equation}\label{coarea}
|Du|(B)=\int_\R P(\{u>t\},B)\, dt,
\end{equation}
with $P_\gamma$ in place of $P$ and $D_\gamma u$ in place of $Du$ in the Gaussian case.

Let us come at a very short discussion of {\em fine properties} of $BV$ functions. 
Observing that, as usual, $BV_{\rm loc}$ functions can be defined as those $L^1_{\rm loc}(\R^d)$ 
functions such that 
\[
V(u,A)=\sup\Big\{\int_{A} u\diver \phi dx:\phi\in C_c^1(A,\R^d),
\|\phi\|_{\infty}\leq 1\Big\}<\infty 
\]
for all bounded open sets $A\subset\R^d$, clearly $BV(\R^d,\gamma)\subset BV_{\rm loc}(\R^d)$,
hence we may confine to $BV_{\rm loc}(\R^d)$ to treat both the Lebesgue and the Gaussian case.  
On the other hand , it is clear that $BV(\R^d)\subset BV(\R^d,\gamma)$ and that in this case 
$D_\gamma u = G_d\, Du$. 

According to the general discussion on approximate limits, we may assume that all 
the functions are approximately continuous in their Lebesgue set, and we may call 
$S_u$ the complement of the Lebesgue set of $u$. Let us list some properties 
of $BV_{\rm loc}$ functions. 

\begin{theorem}\label{fineproperties}
Let $u$ belong to $BV_{\rm loc}(\R^d)$. Then, the following hold:
\begin{itemize}
\item[{\bf (1)}] $S_u$ is an $\LL^d$-negligible and countably $(d-1)$-rectifiable Borel set;
\item[{\bf (2)}] there is $J_u\subset S_u$ such that for every $x\in J_u$ there are 
$u^+(x)\neq u^-(x)\in \R$ and $\nu_u(x)\in {\bf S}^{d-1}$ such that, setting 
\begin{align*}
B_\varrho^+(x)&=B_\varrho(x)\cap\{(y-x)\cdot \nu_u(x)>0\},
\\ 
B_\varrho^-(x)&=B_\varrho(x)\cap\{(y-x)\cdot \nu_u(x)<0\}, 
\end{align*}
the following equalities hold: 
\begin{equation}\label{halflimits}
\begin{array}{l}
\displaystyle{\lim_{\varrho\to 0} \frac{1}{\LL^d(B^+_\varrho)}
\int_{B_\varrho^+(x)} |u(y)-u^+(x)|\, dy = 0,}
\\ \ \\
\displaystyle{\lim_{\varrho\to 0} \frac{1}{\LL^d(B^-_\varrho)}
\int_{B_\varrho^-(x)} |u(y)-u^-(x)|\, dy = 0.} 
\end{array}
\end{equation}
$J_u$ is called {\em approximate jump set}, the values $u^\pm(x)$ {\em approximate one-sided limits}
and $\nu_u(x)$ {\em approximate normal to $J_u$ at $x$}. Moreover, the triple $(u^+(x),u^-(x),\nu_u(x))$ 
is determined up to an exchange between $u^+(x)$ and $u^-(x)$ and a change of sign of $\nu_u(x)$; 
\item[{\bf (3)}] ${\mathcal H}^{d-1}(S_u\setminus J_u)=0$, the functions $x\mapsto u^\pm(x)$, $x\in J_u$, 
are Borel, if $B$ is such that ${\mathcal H}^{d-1}(B)=0$ then $|Du|(B)=0$  and the measure $Du\res J_u$ coincides 
with $(u^+-u^-)\nu_u {\mathcal H}^{d-1}\res J_u$. 
\end{itemize}
\end{theorem}

If $u=\chi_E\in BV_{\rm loc}(\R^d)$ is a characteristic function, we say that the set $E$ has 
{\em locally finite perimeter} and we can say more on the set where the measure $P(E)$ is 
concentrated. Simple examples show that the topological boundary $\partial E$ is too large
(it can be the whole space), hence some suitable relevant subsets should be identified. 
In this connection, the notion of {\em density}, which is slightly weaker than that of approximate 
limit but has a more direct geometric meaning, turns out to be useful. We say that $E\subset\R^d$ 
has density $\alpha\in [0,1]$ at $x\in\R^d$ if
\begin{equation}\label{defdensity}
\lim_{\varrho\to 0}\frac{\LL^d (E\cap B_\varrho(x))}{\LL^d (B_\varrho(x))} = \alpha 
\end{equation}
and in this case we write $x\in E^\alpha$. Of course, if $\aplim_{y\to x}\chi_E(y)=\alpha$ then 
$x\in E^\alpha$. We introduce the {\em essential boundary} 
\[
\partial^*E=\R^d\setminus(E^0\cup E^1)
\] 
and the {\em reduced boundary} ${\mathcal F}E$, defined as follows: $x\in {\mathcal F}E$ if the following 
conditions hold:
\begin{equation}\label{redbdry}
|D\chi_E|(B_\varrho(x))>0\ \forall\, \varrho>0\ {\rm and\ }\exists\ \nu_E(x)=
\lim_{\varrho\to 0} \frac{D\chi_E(B_\varrho(x))}{|D\chi_E|(B_\varrho(x))}
\end{equation}
with $|\nu_E(x)|=1$. If $x\in {\mathcal F}E$, the hyperplane $T(x)=T_{\nu_E(x)}=\{y\in R^d: y\cdot\nu_E(x)=0\}$ 
is the approximate tangent space to ${\mathcal F}E$ as in \eqref{apptgspace}. Indeed,
\begin{equation}\label{apptgredbdry}
\lim_{\varrho\to 0} \frac{E-x}{\varrho} =  \{y\in \R^d: y\cdot\nu_E(x) > 0\}
\end{equation}
locally in measure in $\R^d$. Looking at the properties of $u=\chi_E$, the following inclusions hold:
\[
{\mathcal F}E = J_u \subset E^{1/2} \subset \partial^*E = S_u.  
\]
On the other hand, ${\mathcal H}^{d-1}(\R^d\setminus(E^0\cup E^1\cup E^{1/2}))=0$ 
and in particular ${\mathcal H}^{d-1}(\partial^*E\setminus{\mathcal F}E)=0$.  
For further reference, it is worth noticing that densities are related to the short-time behaviour
of the heat semigroup, i.e., 
\begin{equation}\label{heatdensity}
x\in E^\alpha\quad \Longrightarrow\quad \lim_{t\to 0} W_t\chi_E(x) = \alpha .
\end{equation}
Let us point out now that there are still (at least) two relevant issues 
concerning the infinite dimensional setting, the {\em slicing} and the discussion of {\em embedding 
theorems}, both for Sobolev and $BV$ spaces and the related {\em isoperimetric inequalities}. 
Of course, we are interested here in these arguments in the Gaussian case, and indeed they can be discussed 
directly in the Wiener case, because these results are dimension independent, hence 
there is not a big difference with respect to $(\R^d,\gamma)$ setting.

\section{The Wiener space}\label{secWiener}

In this section we present the measure theoretic and the differential structure which characterize 
the Wiener spaces. After briefly describing the {\em classical Wiener space}, whose elements are 
stochastic processes, we introduce the abstract structure. 

\subsection{Classical Wiener space}\label{classwiener}

For $a\in\R^d$, let $X=C_a([0,1],\R^d)$ be the Banach space of $\R^d$-valued continuous functions 
$\omega$ on $[0,1]$ such that $\omega(0)=a$, endowed with the sup norm and the Borel $\sigma$-algebra 
${\mathcal B}(X)$. Looking at $(X, {\mathcal B}(X))$ as a measurable space, consider the 
{\em canonical process} $B_t(\omega)=\omega(t)$, $0\leq t\leq 1$. Then, there is one probability measure 
$\PP$ (called {\em Wiener measure}) such that $B_t$ is a Brownian motion in $\R^d$ such that $B_0=a$. 
If we want to identify the measure $\PP$, we can exploit the fact that linear and bounded functionals on 
$X$, i.e., Radon measures, can be tought of as random variables. Using the fact that $B_t=\delta_t$ 
and that delta measures are dense in the dual of $X$, it is possible to conclude that $\PP=\NN(a,Q)$ is 
a Gaussian measure with covariance $Q=(q_{hk})$, $q_{hk}=q_h\delta_{hk}$ with  
\[
q_h(\mu,\nu)=\int_0^1\!\int_0^1 s\wedge t\, \mu_h(ds)\, \nu_h(dt),\ \mu,\nu\in{\mathcal M}([0,1],\R^d),
\ h=1,\ldots, d.
\]
Given Borel sets $B_j\in \B(\R^d)$, $j=1,\ldots, m$ and $0=t_0<t_1<\ldots < t_m\leq 1$, define the {\em cylinder} 
\[
C=\{\omega\in X:\ \omega(t_j)\in B_j,\ j=1,\ldots, m\};
\]
we have 
\begin{align}\label{Wienermeas}
\PP(C) = & \int_{B_1} G(t_1,a,x_1)\, dx_1 \int_{B_2} G(t_2-t_1,x_1,x_2)\, dx_2
\\ \nonumber
&\cdots \int_{B_m} G(t_m-t_{m-1},x_{m-1},x_{m})\, dx_m ,
\end{align}
where $G$ is defined in \eqref{heatkernel}. 

For what follows (see \eqref{discreteintbyparts} below), it is important to know for which 
functions $h\in X$ the measure $\PP_h(B)=\PP(h+B)$ is absolutely continuous with respect to 
$\PP$: this happens if and only if $h\in H=X\cap H^1(0,1)$ (Cameron-Martin Theorem \cite{boga}), 
i.e., if and only if $h\in X,\ h'\in L^2(0,1)$. 

As a consequence of the above discussion, the space of the directions which give absolutely 
continuous measures under translation has a natural Hilbert space structure. As we are going 
to see, this is a general fact. 

The same construction of the Wiener measure can be done in the (non separable) space of bounded
Borel functions on $(0,1)$, but by Kolmogorov Theorem (see \cite[Chapter 5]{Wentzell}) the 
Wiener measure concentrates on $C_0([0,1],\R^d)$. 

In this setting, we present a result due to Fukushima, see \cite{fuk99}, which has been 
the starting point of the whole theory, as it highlights a strong connection between the theory of 
perimeters and the stochastic analysis. We use the notation of Section \ref{dirichletforms}.

\begin{theorem}\label{fuk99result}
Given an open set $D\subset \R^d$, the following conditions are equivalent: 
\begin{itemize}
\item[{\bf i)}] $D$ has finite perimeter;
\item[{\bf ii)}] the reflecting Brownian motion $(X_t,{\mathbb P}_x)$ on $\overline D$ is a 
{\em semimartingale}, in the sense that the decomposition  
\[
X_t=X_0+B_t+N_t,
\]
holds, where $B_t$ is the standard $d$-dimensional Brownian motion and each component $N^i_t$ is of 
bounded variation and satisfies the property
\[
\lim_{t\downarrow 0} 
\frac{1}{t}{\mathbb E}\left[\int_0^t \chi_K(X_s) d|N^i_s| \right]<+\infty
\]
for any compact set $K\subset D$.
\end{itemize}
\end{theorem}
The idea is that if $D$ is a set with finite perimeter, then in a weak sense the Brownian motion $B_t$ is
reflected when it reaches the boundary of $D$ since an (approximate) tangent space is defined; using the 
language of processes, the reflecting Brownian motion admits an expression of the form
\[
X_t=X_0+B_t+\int_0^t \nu_D(X_s) dL_s,
\]
where $L_t$ describes the reflection on the boundary; it is the local time, i.e., it is an 
{\em additive functional} with {\em Revuz measure} given by $\cH^{d-1}\res {\mathcal F}D$, that is
\[
\lim_{t\downarrow 0} \frac{1}{t}{\mathbb E}\left[ \int_0^t f(X_s) dL_s \right]
= \int_{{\mathcal F}D} f d\cH^{d-1} 
\]
for continuous $f$. We refer to \cite{Fot94} for the related notions. The idea expressed by this 
theorem is that, since the Brownian motion has trajectories that are not $C^1$ and the 
tangent space to $\partial D$ exists only in an approximate sense, a reflection law 
is not properly defined in terms of classical calculus, but the reflection
properties of the Brownian motion can be described only in a stochastic sense and are contained 
in the additive functional $L_t$, the local time.

Fukushima proves the result for a general $BV$ function $\rho$, by considering the Dirichlet form
\[
\EE(u,v) =\int_{\R^d} \nabla u\cdot \nabla v \rho dx
\]
with associated process $(X_t,{\mathbb P}_x)$. The idea of the proof is to show that the additive functional
\[
A^{[u]}_t=u(X_t)-u(X_0)
\]
admits a semimartingale decomposition
\[
A^{[u]}_t=M^{[u]}_t+N^{[u]}_t,
\]
with $M^{[u]}_t$ a martingale and $N^{[u]}_t$ of bounded variation if and only if
\[
|\EE(u,v)|\leq c\|v\|_\infty,
\]
for some positive constant $c>0$. The particular choice $u(x)=x_i$, the projection 
onto the $i$-th coordinate gives the result.

\subsection{Abstract Wiener spaces}

Let us come to the notion of abstract {\em Wiener space}.  Given a separable Banach space $X$, let 
$\gamma=\NN(0,Q)$ be a nondegenerate centred Gaussian measure on $(X,{\mathcal B}(X))$. As a general
comment, let us point out that a Gaussian measure can be defined in any Banach space, and 
it is always concentrated on a separable subspace, as briefly recalled in the preceding 
subsection. Moreover, a consequence of Fernique's theorem, see \eqref{fernique}, is that any 
$x^*\in X^*$ defines a function $x\mapsto \scal{x}{x^*}$ belonging to $L^p(X,\gamma)$ for 
all $p\geq 1$. In particular, we may think of any $x^*\in X^*$ as an element of $L^2(X,\gamma)$. 
Let us denote by $R^*:X^*\to L^2(X,\gamma)$ the embedding, $R^*x^*(x)=\scal{x}{x^*}$. The closure 
of the image of $X^*$ in $L^2(X,\gamma)$ under $R^*$ is denoted $\mathscr{H}$ and is called the 
{\em reproducing kernel} of the Gaussian measure $\gamma$. The above definition is motivated by 
the fact that if we consider the operator $R:\mathscr{H}\to X$ whose adjoint is $R^*$, then 
\begin{equation}\label{i}
R\hat h=\int_X \hat h(x)xd\gamma(x),\qquad \hat{h}\in {\mathscr H}
\end{equation}
(Bochner integral). In fact, denoting by $[\cdot,\cdot]_{\mathscr{H}}$ the inner product in $\mathscr{H}$
and by $|\cdot|_{\mathscr{H}}$ the norm, the equality
\[
[\hat{h},R^*x^*]_{\mathscr H}=\int_X\hat{h}(x)\scal{x}{x^*}d\gamma(x)
=\Big\langle\int_X\hat{h}(x)xd\gamma(x),x^*\Big\rangle ,
\]
that holds for all $x^*\in X^*$, implies \eqref{i}.  With the definition of $R,\,R^*$ we obtain 
directly by \eqref{CovOp} the decomposition $Q=RR^*$:
\[
\scal{RR^*x^*}{y^*}=[R^*x^*,R^*y^*]_{\mathscr H}
=\int_X \scal{x}{x^*}\scal{x}{y^*}d\gamma(x)
=\scal{Qx^*}{y^*}.
\] 
The space $H=R\mathscr{H}$ is called the {\em Cameron-Martin space}; it is a Hilbert space, 
dense in $X$ because $\gamma$ is nondegenerate, with inner product defined by
\[
[h_1,h_2]_{H}=[\hat h_1,\hat h_2]_{\mathscr{H}}
\]
for all $h_1,h_2\in H$, where $h_i=R\hat h_i$, $i=1,2$, and norm $|\cdot|_H$. As recalled in Subsection 
\ref{Measuretheory}, $Q$ is a compact operator. The same holds for $R$ and $R^*$, hence the 
embeddings $X^*\hookrightarrow {\mathscr H}$, $H\hookrightarrow X$ are compact.   
Given the elements $x^*_1,\ldots,x^*_m$ in $X^*$, we denote by $\pi_{x^*_1,\ldots,x^*_m}:X\to \R^m$
the finite dimensional projection of $X$ onto $\R^m$ induced by the elements $x^*_1,\ldots,x^*_m$, 
that is the map
\[
\pi_{x^*_1,\ldots,x^*_m}x=(\scal{x}{x^*_1},\ldots,\scal{x}{x^*_m}),
\]
also denoted by $\pi_m:X\to \R^m$ if it is not necessary to specify the elements 
$x^*_1,\ldots,x^*_m$. The symbol $\FF C_b^k(X)$ denotes the space of $k$ times continuously 
differentiable {\em cylindrical} functions with bounded derivatives up to the order $k$, that is: 
$u\in \FF C_b^k(X)$ if there are $m\in\N$, $x^*_1,\ldots, x^*_m\in X^*$ and $v\in C_b^k(\R^m)$ 
such that $u(x)=v(\pi_mx)$. We denote by ${\mathcal E}(X)$ the cylindrical $\sigma$-algebra 
generated by $X^*$, that is the $\sigma$-algebra generated by the sets of the form $E=\pi_m^{-1}B$ 
with $B\in {\mathcal B}(\R^m)$. Since $X$ is separable, ${\mathcal E}(X)$ and $\B(X)$ coincide, see 
\cite[Theorem I.2.2]{VTC}, even if we fix a sequence $(x^*_j)\subset X^*$ which separates the points 
in $X$ and use only elements from that sequence to generate $\pi_m$. We shall make later on some 
special choices of $(x^*_j)$, induced by the Gaussian probability measure $\gamma$ in $X$.
Using the embedding $R^*X^*\subset\mathscr{H}$, we say that a family 
$\{x^*_j\}$ of elements of $X^*$ is orthonormal if the corresponding family $\{R^*x^*_j\}$ is 
orthonormal in $\mathscr{H}$. It can be proved that $\gamma(H)=0$, see \cite[Theorem 2.4.7]{boga}
Since $X$ and $X^*$ are 
separable, starting from a sequence in $X^*$ dense in $H$, we may construct 
an orthonormal basis $(h_j)$ in $H$ with $h_j=Qx^*_j$. Set also $H_m={\rm span}\{h_1,\ldots,h_m\}$,
and define $X^\perp ={\rm Ker}\, \pi_{x^*_1,\ldots,x^*_m}$ and $X_m$ the ($m$-dimensional) 
complementary space. Accordingly, we have the canonical decomposition $\gamma=\gamma_m\otimes\gamma^\perp$ 
of the measure $\gamma$; notice also that these Gaussian measures are {\em rotation invariant}, i.e.,
if $\varrho:X\times X\to X\times X$ is given by 
$\varrho(x,y)=(\cos\vartheta x+\sin\vartheta y,-\sin\vartheta x + \cos\vartheta y)$ 
for some $\vartheta\in \R$, then $\varrho_\#(\gamma\otimes\gamma)=\gamma\otimes\gamma$ 
and the following equality holds:
\begin{equation} \label{rotinv}
\int_{X}\int_{X} u(\cos\vartheta x + \sin\vartheta y)d\gamma (x)d\gamma (y)
=\int_{X} u(x)d\gamma (x),
\end{equation}
$u\in L^1(X,\gamma)$, which is obtained by the above relation by integrating the function 
$u\otimes \one$ on $X\times X$. Notice that if $X$ is decomposed as $X=X_m\oplus X^\perp$, the 
same formula holds in $X_m$ and $X^\perp$ separately, with measures $\gamma_m$ and $\gamma^\perp$.  

For every function $u\in L^1(X,\gamma)$, if $\{h_j\}$ is 
an orthonormal basis of $H$, its {\em canonical cylindrical approximations} $u^m$ are defined
as the conditional expectations relative to the $\sigma$-algebras $\F_m=\pi_m^{-1}(\B (\R^m))$,   
\begin{equation}\label{cancyl}
u^m=\E (u|\F_m) = \E_m u\quad {\rm s.t.}\quad \int_Aud\gamma=\int_Au^md\gamma
\end{equation}
for all $A\in\F_m$. Then, $u^m\to u$ in $L^1(X,\gamma)$ and $\gamma$-a.e. 
(see e.g. \cite[Corollary 3.5.2]{boga}). More explicitly, we set
\[
\E_m u(x)=\int_{X}u(P_mx+(I-P_m)y)d\gamma(y)
=\int_{X^\perp}u(P_mx+y')d\gamma^\perp(y'),
\] 
where $P_m$ is the projection onto $X_m$. Notice that the restriction of $\gamma$ to $\F_m$ is 
invariant under translations along all the vectors in $X^\perp$, hence we may write 
$\E_mu(x)=v(P_mx)$ for some function $v\in L^1(X_m,\gamma_m)$, and, with 
an abuse of notation, $\E_mu(x_m)$ instead of $\E_mu(x)$. 
 
The importance of the Cameron-Martin space relies mainly on the fact that the translated measure
\[
\gamma_h(B)=\gamma(B-h),\qquad B\in \B(X),\quad  h\in X
\]
is absolutely continuous with respect to $\gamma$ if and only if $h\in H$ and in this case, with the usual 
notation $h=R\hat h$, $\hat h\in{\mathscr H}$, we have, see e.g. \cite[Corollary 2.4.3]{boga},
\begin{equation} \label{cameronmartin}
d\gamma_h(x)=\exp\Big\{\hat h(x)-\frac{1}{2}|h|_H^2\Big\}d\gamma(x).
\end{equation}
Let us look for the basic integration by parts formula in the present context, that generalises 
\eqref{gaussintbyparts} and allows to define weak derivatives and $BV$ functions.  
For $h\in X$, define 
\[
\partial_hf(x)=\lim_{t\to 0}\frac{f(x+th)-f(x)}{t}
\]
(whenever the limit exists); we look for an operator $\partial^*_h$ such that 
for every $f,g\in \FF C_b^1(X)$ the equality 
\begin{equation}\label{intbyparts1}
\int_X g(x) \partial_h f(x)d\gamma(x)=-\int_X f(x)\partial^*_h g(x)d\gamma(x)
\end{equation}
holds. Starting from the incremental ratio, we get
\begin{align}\label{discreteintbyparts}
\int_X \frac{f(x+th)-f(x)}{t} g(x) d\gamma(x)=&
-\int_X f(y) \frac{g(y)-g(y-th)}{t} d\gamma_{th}(y)
\\ \nonumber
&+\int_X f(x)g(x)d\mu_t(x)
\end{align}
where $\mu_t=\frac{1}{t}\Bigl(\NN(th,Q)-\NN(0,Q)\Bigr).$ 
{}From the Cameron-Martin formula \eqref{cameronmartin} we know that $\mu_t\ll\gamma$
if and only if $h\in H$. In this case, we can use \eqref{cameronmartin}
and pass to the limit by dominated convergence as $t\to 0$, getting \eqref{intbyparts1}
with
\[ 
\partial^*_hg(x)=\partial_hg(x)-g(x)\hat{h}(x), 
\]
where as usual $h=R\hat{h}$. Such notions can be extended to the more general class of
{\em differentiable measures}, see \cite{bogadm}. Let us now define the gradient and the divergence 
operators. For $f\in \FF C^1_b(X)$, the $H$-gradient of $f$, denoted by 
$\nabla_H f$, is the map from $X$ into $H$ defined by
\[
[\nabla_H f(x),h]_{H}= \partial_h f(x),
\quad h\in H,
\]
where $\partial_h f(x)$ is defined as before.
Notice that if $f(x)=f_m(\pi_mx)$ with $f_m \in C^1(\R^m)$, then
\[
\partial _h f(x)=\nabla f_m(\pi_mx) \cdot \pi_m h.
\]
If we fix an orthonormal basis $\{h_j\}_{j\in \N}$ of $H$, we can write
\[
\nabla_Hf(x)=\sum_{j\in \N}\partial_j f(x)h_j, \qquad \partial _j=\partial_{h_j},
\]
where it is important to notice that the directional derivative $\partial_h$ is computed
by normalising $h$ with respect to the norm in $H$. 
Considering the space $\FF C^1_b(X,H)$, 
we may define $-\diver_H$, the adjoint operator of $\nabla_H$, as the
linear map from $\FF C^1_b(X,H)$ to $\FF C_b(X)$ such that
\[
\diver_H\phi(x)=\sum_{j\in \N} \partial^*_j\phi_j(x)
=\sum_{j\in\N}\partial_j\phi_j(x)-\phi_j(x)\hat{h}_j(x),\quad
\phi_j=[\phi,h_j]_H.
\] 

\subsection{Hausdorff measures}\label{SubsHausdorff}

The definition of Hausdorff measures in Wiener spaces goes back to \cite{feydelpra} and is based 
on a finite dimensional approximation. If $F\subset X$ is an $m$-dimensional subspace of $H$, 
$B\subset F$, recall that we are denoting by ${\mathcal S}^{k}(B)$ the {\em spherical $k$-dimensional 
Hausdorff measure} of $B$. We stress that the balls used in the minimisation above are understood 
with respect to the $H$ distance and we do not emphasise the dependence on $F$. Occasionally we 
canonically identify $F$ with $\R^m$, choosing a suitable orthonormal basis.

Let $F\subset QX^*$ be an $m$-dimensional subspace of $H$. We
denote by $z=\pi_F(x)$ the canonical projection induced by an
orthonormal basis $e_j=Qe_j^*$ of $F$, namely
\[
\pi_F(x)=\sum_{j=1}^m\langle e_j^*,x\rangle e_j
\]
and set $x=y+z$, so that $y=x-\pi_F(x)$ belongs to ${\rm Ker}(\pi_F)$, the kernel of $\pi_F$. 
This decomposition induces the factorization $\gamma=\gamma^\perp\otimes\gamma_F$ with $\gamma_F$
standard Gaussian in $F$ and $\gamma^\perp$ Gaussian in ${\rm Ker}(\pi_F)$ (whose Cameron-Martin 
space is $F^\perp$).

Following \cite{feydelpra}, we can now define spherical $(\infty-1)$-dimensional Hausdorff measures 
in $X$ relative to $F$ by
\begin{equation}\label{pre1}
{\mathcal S}^{\infty-1}_F(B)= 
\int_{{\rm Ker}(\pi_F)}^*\int_{B_y}G_m(z)\,d{\mathcal S}^{m-1}(z) \,d\gamma^\perp(y)
\qquad \forall B\subset X. 
\end{equation}
Here, for $y\in {\rm Ker}(\pi_F)$, by $B_y$ we denote the {\em section} or {\em slice} 
\begin{equation}\label{pausa}
B_y=\left\{z\in F:\ y+z\in B\right\}.
\end{equation}
The internal integral in \eqref{pre1} is understood in the Choquet sense, namely
\[
\int_{B_y}G_m(z)\,d{\mathcal S}^{m-1}(z)=
\int_0^\infty{\mathcal S}^{m-1}\bigl(\{z\in B_y:\ G_m(z)>\tau\}\bigr)\,d\tau.
\]
If $B_y\in{\mathcal B}(F)$, as it happens in the case $B\in{\mathcal B}(X)$, the integral 
reduces to a standard one. Furthermore, we have used the outer integral in order to avoid 
the issue of the measurability of the map $y\mapsto\int_{B_y}G_m\,d{\mathcal S}^{m-1}$.
The next basic additivity result is proved in \cite{feydelpra}. 

\begin{lemma}\label{sigmadd}
${\mathcal S}^{\infty-1}_F$ is a $\sigma$-additive Borel measure on
${\mathcal B}(X)$. In addition, for all Borel sets $B$ the map
$y\mapsto \int_{B_y}G_m\,d{\mathcal S}^{m-1}$ is $\gamma^\perp$-measurable
in ${\rm Ker}(\pi_F)$.
\end{lemma}

A remarkable fact is the monotonicity of ${\mathcal S}^{\infty-1}_F$ with respect to $F$, 
which crucially depends on the fact that we are considering spherical Hausdorff measures.

\begin{lemma} \label{lmono}
${\mathcal S}^{\infty-1}_F\leq{\mathcal S}^{\infty-1}_G$ on ${\mathcal B}(X)$ whenever
$F\subset G$.
\end{lemma}

The above property has been pointed out in \cite{feydelpra}, relying on \cite[2.10.27]{federer}.
We refer to \cite[Lemma 3.1]{AMP} for details. It follows from Lemma~\ref{lmono} 
that the following definition of {\em spherical $(\infty-1)$-Hausdorff measure} 
${\mathcal S}^{\infty-1}$ in ${\mathcal B}(X)$ is well-posed; we set 
\begin{equation}\label{Sinfty-1}
{\mathcal S}^{\infty-1}(B)=\sup_F{\mathcal S}^{\infty-1}_F(B)=\lim_F{\mathcal S}^{\infty-1}_F(B),
\end{equation}
the limits being understood in the directed set of finite-dimensional subspaces of $QX^*$. 
A direct consequence of Lemma \ref{sigmadd} is that ${\mathcal S}^{\infty-1}$ is $\sigma$-additive 
on ${\mathcal B}(X)$. This measure does not coincide with the one of \cite{feydelpra}, since we are 
considering only subspaces in $H$ generated by elements of $QX^*$. Our approach is a bit simpler 
because the corresponding projections are continuous, whereas general orthogonal decompositions 
of $H$ give merely measurable projections, so that some technical points 
related to removing sets of small capacity has to be addressed. 

\subsection{Sobolev spaces and isoperimetric inequality}\label{isoper_ineq}

There are several possible definitions of Sobolev spaces on Wiener spaces. 
Since the operator $\nabla_H$ is a closable operator in 
$L^p(X,\gamma)$, one may define the Sobolev space ${\mathbb D}^{1,p}(X,\gamma)$ 
as the domain of the closure of $\nabla_H$ in $L^p(X,\gamma)$\footnote{Notice that the space 
denoted by ${\mathbb D}^{1,p}(X,\gamma)$ by Fukushima is denoted by $W^{p,1}(X,\gamma)$ 
in \cite{boga}.}. Another possible definition, which is closer to our point of view, is 
based on the integration by parts formula \eqref{intbyparts1}: $f\in L^p(X,\gamma)$ is in
$W^{1,p}(X,\gamma)$ if there is $F\in L^p(X,\gamma;H)$ such that \eqref{intbyparts1} 
holds with $[F,h]_H$ in place of $\partial_hf$ and any $g\in \FF C^1_b(X,H)$. In this case,
we denote $F$ by $\nabla_Hf$. Anyway, the spaces $W^{1,p}$ and ${\mathbb D}^{1,p}$ coincide, 
see \cite[Section 5.2]{boga}. This approach requires some further explanations in the case
$p=1$, as we shall see at the end of this subsection. 

The {\em Gaussian isoperimetric inequality} says the following, see \cite{ledoux}. Let $E\subset X$, 
and set $B_r=\{x\in H:\|x\|_{H}<r\},\ E_r=E+B_r$;  then
\[
\Phi^{-1}(\gamma(E_r))\geq \Phi^{-1}(\gamma(E))+r, \quad {\rm where}\ 
\Phi(t):=\int_{-\infty}^t\frac{e^{-s^2/2}}{\sqrt{2\pi}}ds.
\] 
We sketch here why this inequality implies the isoperimetric inequality. We
introduce the function
\[
\mathscr{U}(t):=(\Phi'\circ\Phi^{-1})(t)
\approx t\sqrt{2\log(1/t)},\ t\to 0.
\] 
Since $\mathscr{U}(t)=\mathscr{U}(1-t)$, the function $\mathscr{U}$ has the same behaviour 
as $t\to 1$, $\mathscr{U}(t)\approx (1-t)\sqrt{2\log(1-t)}$. Notice that $\Phi(t)$ is the 
volume of the halfspace $\{\hat{h}(x)<t\}$ and that $\cU(t)$ is the perimeter of a halfspace 
of volume $t$.

{}From the above estimate for $\Phi^{-1}(\gamma(E_r))$ we obtain that
\begin{align*} 
\gamma(E_r)&\geq \Phi(\Phi^{-1}(\gamma(E))+r) = \gamma(E)+r \Phi'(\Phi^{-1}(\gamma(E)))+o(r) 
\\
&= \gamma(E)+r\mathscr{U}(\gamma(E))+o(r),
\end{align*}
and then
\[
\liminf_{r\to 0} \frac{\gamma(E_r)-\gamma(E)}{r} \geq \mathscr{U}(\gamma(E)).
\] 
The quantity on the left hand side is related to the Minkowski content of the set $E$ constructed 
using the Cameron-Martin balls, although negligible. For instance, if $X=\R^d$, $\gamma=G_d\LL^d$ the 
standard centred Gaussian measure on $\R^d$ and $E$ a set with smooth boundary, then
\[ 
P_{\gamma}(E)=\lim_{r\to 0}\frac{\gamma(E_r)-\gamma(E)}{r} 
\geq \mathscr{U}(\gamma(E)).
\]  
It is also possible to prove in this case that equality holds if $E$ is a hyperplane;
this skecth of the isoperimetry property of hyperplanes is essentially the proof contained
in \cite{ledoux}. The original proof of the isoperimeric properties of hyperplanes in the 
finite dimensional Gaussian space has been estabilished first in \cite{sudtsi};
since the isoperimetric function does not depend on the space dimension, the
same proof can be extended to the infinite dimensional case. In \cite{Ehr83Sym}, again
in the finite dimensional case, it is proved that hyperplanes are isoperimetric by using a
symmetrisation argument; also in this case, the proof implies that hyperplanes
are isoperimetric in the infinite dimensional case. The proof  that hyperplanes are the 
{\em unique} isoperimetric sets is rather recent and is contained in \cite{CarKer01OnT}.
Let us also point out that the right Minkowski content uses enlargements $E_r$ of the set 
$E$ with respect to balls of $H$ and not of $X$. The reason of this can be explained as 
follows: the Gaussian measure $\gamma$ introduces an anisotropy on $X$ due to the covariance 
operator $Q$. This anisotropy is compensated in the definition of total variation and 
perimeter by the gradient $\nabla_H$, since it is defined using vectors that have unit 
$H$-norm. The corresponding compensation in the computation of the Minkowski content is 
achieved by using the balls of $H$.

The isoperimetric inequality implies also the following:
\[ 
\|\nabla_H f\|_{L^1}\geq \int_0^{\infty}\mathscr{U}\left(\gamma(\{|f|>s\})\right)ds,
\] 
and it follows that if $\nabla_Hf\in L^1(X,\gamma)$ then $u$ belongs to the {\em Orlicz space}
\begin{equation}\label{defLLogL}  
\LLog=\{u:X\to \R :\ A_{1/2}(|u|)\in L^1(X,\gamma)\}, 
\end{equation}
where $A_{1/2}(t)=\int_0^t\log^{1/2}(1+s)ds$. This is important in connection to the 
integration by parts formula \eqref{intbyparts1}, because for general $f\in L^1(X,\gamma)$ 
the product $\hat{h}fg$ is not summable. But, thanks to Fernique theorem, the linear function 
$\hat{h}$ belongs to the Orlicz space defined through the complementary $N$-function of 
$A_{1/2}$, $\psi(t)=\int_0^t(e^{s^2}-1)ds$, i.e., $\psi(\lambda|h|_H)<\infty$ for some $\lambda>0$. 
As a consequence, if $\nabla_Hf\in L^1(X,\gamma)$ then $f\in\LLog$, the product $\hat{h}fg$ is 
summable, \eqref{intbyparts1} does make sense and the embedding of ${\mathbb D}^{1,1}(X,\gamma)$ 
into $\LLog$ follows, see \cite[Proposition 3.2]{fuk2000_2}. 
 
\subsection{The Ornstein-Uhlenbeck semigroup}\label{OUsection}

Let us consider the Ornstein-Uh\-len\-beck semigroup $(T_t)_{t\geq 0}$, defined 
pointwise by Mehler's formula, which generalises \eqref{OUsgrpinRd}:  
\begin{equation}\label{OUfor}
T_t u(x)=\int_X u\left(e^{-t}x +\sqrt{1-e^{-2t}}y\right)d\gamma(y)
\end{equation}
for all $u\in L^1(X,\gamma)$, $t>0$. Unlike the heat semigroup, the Ornstein-Uhlenbeck semigroup  
$T_t$ does not map $L^1(X,\gamma)$ into ${\mathbb D}^{1,1}(X,\gamma)$. But, $T_t$ is strongly 
continuous in $\LLog$ and it follows from \eqref{rotinv} that $T_t u\in  {\mathbb D}^{1,1}(X,\gamma)$ 
for any $u\in \LLog$, see \cite[Proposition 3.6]{fuk2000_2}. Moreover, it is a contraction semigroup 
in $L^p(X,\gamma)$ for every $p\in [1,+\infty]$ (and hence also in $\LLog$) and self-adjoint in 
$L^2(X,\gamma)$. Moreover, the following commutation relation holds for any $u\in {\mathbb D}^{1,1}(X,\gamma)$
\begin{equation}\label{commutation}
\nabla_H T_t u=e^{-t}T_t \nabla_H u,\qquad t>0.
\end{equation}
Therefore, we get 
\[ 
\nabla_HT_{t+s}u=\nabla_HT_t (T_s u)
= e^{-t}T_t \nabla_H T_s u,
\] 
for any $u\in \LLog$, see \cite[Proposition 5.4.8]{boga}. 
It also follows from \eqref{commutation} that
\begin{equation}\label{commdiver}
\int_X T_tf\diver_H\phi d\gamma=e^{-t}\int_Xf\diver_H(T_t\phi) d\gamma,
\end{equation}
for all $f\in L^1(X,\gamma),\ \phi\in {\mathscr F}C^1_{b}(X,H)$, see \cite{AMMPPhysicad}. 
Another important consequence of \eqref{commutation} is that if 
$u\in {\mathbb D}^{1,1}(X,\gamma)$ then 
\begin{equation}\label{stronglySobolev}
\lim_{t\to 0}\|\nabla_HT_tu - \nabla_Hu \|_{L^1(X,\gamma)}=0.
\end{equation}
Finally, notice that if $u^m$ are the canonical cylindrical approximations
of a function $u\in \LLog$ defined in \eqref{cancyl} then the following 
inequality holds, see e.g. \cite{AMMPPhysicad}
\begin{equation}\label{monotonicity}
\int_X|\nabla_HT_tu^m|_Hd\gamma \leq \int_X|\nabla_HT_tu|_Hd\gamma
\quad \forall\, t>0.
\end{equation}
We end this brief discussion on the Ornstein-Uhlenbeck semigroup by presenting the related 
Ornstein-Uhlenbeck {\em process} in the Wiener space. Of course, this is close to the 
finite dimensional case, with important modifications. First, we define the {\em cylindrical 
Brownian motion} in $X$ as an $X$-valued continuous process $B^H_t$ such that for every 
$x^*\in X^*$ with $|Qx^*|_H=1$ the one-dimensional process $\langle x^*,B^H_t\rangle$ is 
Wiener. After extending the notion of stochastic integral to the case of a cylindrical Brownian 
motion, we may deal with SDEs in $X$.  The Ornstein-Uhlenbeck process is given by 
\[
\xi_t=e^{-t/2}\xi_0 + \int_0^t e^{(s-t)/2}\, dB^H_t
\]
and, as in $\R^d$, it is the solution of the Cauchy problem for the Langevin equation 
\[
d\xi_t=-\frac{1}{2}\xi_t\, dt+dB^H_t,\qquad \xi_0\ \ {\rm given\ r.v.}, 
\]
where $B^H_t$ is a cylindrical Brownian motion. If the law of $\xi_0$ is $\delta_x$ for $x\in X$,
denoting by $\xi_t^x$ the corresponding solution, we have the usual equality 
$T_tf(x)=\E[f(\xi_t^x)]$. 

\section{$BV$ functions in the Wiener space}\label{secBVWiener}

A definition of $BV$ functions in abstract Wiener spaces has been given by M. Fukushima 
in \cite{fuk2000_1}, M. Fukushima and M. Hino in \cite{fuk2000_2}, and is based upon Dirichlet form 
theory quoted in Subsection \ref{dirichletforms}. In \cite{AMMPPhysicad}, \cite{AMMP} the main aim has 
been to compare the finite and infinite dimensional theory of $BV$ functions from a purely analytical 
point of view, closer to the classical setting. After collecting, in the preceding section, the tools 
we need, we pass now to the definition of $BV$ functions in the abstract Wiener space setting.
We denote by $\cM(X,H)$ the space of all $H$-valued finite measures $\mu$ on $\B(X)$.

\begin{defn}\label{defBVWiener}
Let $u\in \LLog$. We say that $u$ has boun\-ded variation in $X$ and we set $u\in BV(X,\gamma)$ 
if there exists $\mu \in \cM(X,H)$
such that for any $\phi\in {\mathscr F}C^1_b(X,H)$ we have
\begin{equation}\label{BVintbyparts}
\int_X u(x) \partial^*_j\phi(x) d\gamma(x)=-\int_X \phi(x) d\mu_j(x)\qquad \forall j\in \N,
\end{equation}
where $\mu_j=[h_j,\mu]_H$. In particular, if $u=\chi_E$ and $u\in BV(X,\gamma)$, then we say
that $E$ has finite perimeter.
\end{defn}
Notice that, as in the Sobolev case ${\mathbb D}^{1,1}(X,\gamma)$, the assumption $u\in \LLog$ 
gives a meaning to \eqref{intbyparts1}, as discussed in Subsection \ref{isoper_ineq}. 
Moreover, in the previous definition we have required that the measure $\mu$ is defined on the 
whole of ${\B}(X)$ and is $\sigma$-additive there. Since cylindrical functions generate the Borel 
$\sigma$-algebra, the measure $\mu$ verifying \eqref{BVintbyparts} is unique, and will be denoted
$D_\gamma u$ as in the finite dimensional Gaussian case. The total variation measure is denoted 
as usual by $|D_\gamma u|$. We also let $P_\gga(E):=|D_\gamma \chi_E|(X)$ be the (Gaussian) perimeter 
of a subset $E$ of $X$ and we set, as in the finite dimensional case, $P_\gamma(E,\cdot)=|D\chi_{E}|(\cdot)$. 

We state now a characterization of $BV(X,\gamma)$ functions analogous to Theorem \ref{findimeq} and
the discussion which follows. 

\begin{theorem}\label{BVinfinitodim}
Given $u\in \LLog$, the following are equivalent:
\begin{itemize}
\item[{\bf (1)}] $u$ belongs to $BV(X,\gamma)$;
\item[{\bf (2)}] the quantity 
\begin{equation*}
V_\gamma(u):=\sup\Bigl\{\int_X u\diver_H \Phi\, d\gamma;\,\Phi\in {\mathscr F}C^1_b(X,H),
|\Phi(x)|_H\leq 1\ \forall\,x\in X \Bigr\} 
\end{equation*}
is finite;
\item[{\bf (3)}] the quantity 
\begin{equation*}
L_\gamma(u)=\inf\Big\{\liminf\limits_{n\to \infty}\int_X
|\nabla_Hu_n|_Hd\gamma: u_n\in{\mathbb D}^{1,1}(X,\gamma),
u_n\stackrel{L^1}{\rightarrow}u\Big\}
\end{equation*}
is finite;
\item[{\bf (4)}] the quantity
\begin{equation}\label{degiorgilim}
{\mathcal T}[u]=\lim_{t\downarrow 0}\int_X |\nabla_H T_tu|_H d\gamma
\end{equation}
is finite.
\end{itemize}
Moreover, $|D_\gamma u|(X)=V_\gamma(u)=L_\gamma(u)={\mathcal T}[u]$. 
\end{theorem}
As in the finite dimensional case, see \eqref{BVconDirform}, $u\in BV(X,\gamma)$ if and
only if there is $C>0$ such that 
\[
\Big|\int_X [\nabla_H\Phi,h]_H\ u\ d\gamma \Big| \leq C \|\Phi\|_\infty
\]
for all $\Phi\in \FF C^1_b(X,H)$. 
The proof of Theorem \ref{BVinfinitodim} is contained in \cite{fuk2000_1}, \cite{fuk2000_2}, and also 
in \cite{AMMP}. The proof in the latter reference relies on a {\em slicing} argument, 
a technique that has proved to be very useful in the finite dimensional case and we shall
use later. For $\nu\in\bigcup_m H_m$, denote by $\partial_\nu$ and $\partial^*_\nu$ the 
differentiation operator and its adjoint, respectively, and the {\em directional total 
variation along $\nu$} as
\begin{equation}\label{gaussianRm1}
V^\nu_\gamma(u)=\sup\Big\{\int_X u\partial^*_\nu \phi d\gamma:
\phi\in {\mathscr F}^\nu C_c^1(X),\ |\phi(x)|\leq 1
\ \forall\, x\in X\Big\},
\end{equation}
where $\phi\in {\mathscr F}^\nu C_c^1(X)$ means that $\phi(x)=v(\langle x, x^*\rangle)$ with 
$v\in C^1_c(\R)$ and $\nu=Qx^*$.
Riesz theorem shows that $V^\nu_\gamma(u)$ is finite if and only if the integration by parts formula
\begin{equation}\label{amb1}
\int_X u\partial^*_\nu \phi d\gamma=-\int_X\phi d\mu_\nu \qquad\forall
\phi\in C_b^1(X)
\end{equation}
holds for some real-valued measure $\mu_\nu$ with finite total variation, that we denote by 
$D^\nu_\gamma u$; if this happens, $|\mu_\nu|(X)$ coincides with $V^\nu_\gamma(u)$. Finally,
\begin{equation}\label{vienna2}
V^\nu_\gamma(u)=\lim_{m\to\infty} V^\nu_\gamma({\mathbb E}_mu).
\end{equation}
Once a direction $\nu=Qx^*\in H$ is fixed, let $\pi_{\nu}(x)=\langle x,x^*\rangle$ be the induced 
projection and let us write $x\in X$ as $y+\pi_{\nu}(x)\nu$. Then, denoting by $K$ the kernel of
$\pi_{\nu}$, $\gamma$ admits a product decomposition $\gamma=\gamma^\perp\otimes\gamma_1$ with 
$\gamma^\perp$ Gaussian in $K$. For $u:X\to\R$ and $y\in K$ we define the function
$u_y:\R\to\R$ by $u_y(t)=u(y+t\nu)$. The following slicing theorem holds 

\begin{theorem}\label{thmSlicing}
Let $u\in\LLog$ and let $\nu\in \bigcup_mH_m$; then
\[
V^\nu_\gamma(u)=\int_K V_{\gamma_1}(u_y)\,d\gamma^\perp(y).
\]
In particular, the directional total variation of $u$ is independent
of the choice of the basis and makes sense for all $h\in H$.
\end{theorem}

The coarea formula \eqref{coarea} holds as well in Wiener spaces and can be proved by 
following {\em verbatim} the proof of \cite[Section 5.5]{evgar}: if $u\in BV(X,\gamma)$, 
then for a.e. $t\in\R$ the level set $\{u>t\}$ has finite perimeter and for every Borel 
set $B\subset X$ the following equality holds:
\begin{equation}\label{coareaform}
|D_\gamma u|(B)=\int_\R P_\gamma(\{u>t\},B) dt.
\end{equation}

We end this section with a recent example of application in the classical Wiener space, 
see \cite{trevisan}.
 
\begin{example}\label{extrevisan}
{\rm Let us fix a time $t\in [0,1]$ and consider the classical Wiener space $X=C_0([0,1],\R)$, 
see Subsection \ref{classwiener}. Define 
\[
M_t=\sup\{ B_s,\ 0\leq s\leq t\}.
\] 
It is well-known that $M_t\in {\mathbb D}^{1,p}(X,{\mathbb P})$, but $\nabla_HM_t$ is not differentiable. 
Nevertheless, $\nabla_HM_t$ belongs to $BV(X,{\mathbb P})$, i.e., there exists a $H\otimes H$-valued 
measure $D^2M_t$ such that 
\[
\int_X [\Phi h_1\otimes h_2, D^2M_t]_{H\otimes H} = 
\int_X M_t\partial^*_{h_1}\partial^*_{h_2}\Phi\, d\gamma
\]
for every $\Phi\in \FF C^2_b(X)$, $h_1,h_2\in H$. 
Moreover, the measure $|D_{\mathbb P}\nabla_HM_t|$ is concentrated on the trajectories that attain their
maximum exactly twice, hence, in particular, all these measures are singular with respect to ${\mathbb P}$.
}\end{example}

\section{Fine properties of sets with finite perimeter}\label{secfineprop}

We show in this section how is it possible to generalise in the infinite-dimensional setting 
the properties listed in Theorem~\ref{fineproperties}; we restrict our attention 
to the case of sets with finite perimeter, so that we can use the geometric meaning
of points of density stated by formula~\eqref{defdensity} to give a suitable notion of
boundary of a set.

It is worth noticing that in the infinite-dimensional setting things do not work as 
well as for the Euclidean case; Preiss~\cite{Preiss} gave an example of an infinite-dimensional
Hilbert space $X$, a Gaussian measure $\gamma$ and a set $E\subset X$ such that
$0<\gamma(E)<1$ and
\begin{equation}\label{density1}
\lim_{\varrho\to 0}\frac{\gamma(E\cap B_\varrho(x))}{\gamma(B_\varrho(x))}=1,
\qquad \forall x\in X.
\end{equation}
In the same work, it is also shown that if the eigenvalues of the covariance $Q$ decay
to zero sufficiently fast, then it is possible to talk about density points; in some sense, 
the requirement on the decay gives properties of $X$ closer to the finite-dimensional case.
For these reasons, in general the notion of point of density as given in \eqref{density1} 
is not a good notion. 

In the infinite-dimensional setting, the idea is to use the factorization $\gamma=\gamma^\perp \otimes \gamma_F$,
for $F\subset QX^*$ an $m$-dimensional space, described in Subsection~\ref{SubsHausdorff}.

\begin{defn}[Essential boundary relative to $F$]
If we write $X=F\oplus {\rm Ker}(\pi_F)$, we recall by \eqref{pausa} the definition of the slice 
of $E$ in direction $F$
\[
E_y =\{ z\in F: y+z\in E\}\subset F;
\]
the {\emph essential boundary of $E$ relative to $F$} is then defined as
\[
\partial^*_F E=\{ x=y+z: z\in \partial^* (E_y)\}.
\]
\end{defn}
It is not difficult to show that $\partial_F^*E$ is a Borel set; moreover, in order to pass 
from the finite dimensional space $F$ to the whole of the Cameron-Martin space $H$, we need 
the following property. 

\begin{lemma}\label{finitedim}
Let $G\subset QX^*$ be a $k$-dimensional Hilbert space, let $F\subset G$ be an 
$m$-dimensional subspace and let $E$ be a set with finite perimeter
in $G$. Then, with the orthogonal decomposition $G=F\oplus L$ and
the notation
$$
E_w:=\left\{z\in F:\ w+z\in E\right\}\qquad\qquad w\in L,
$$
we have that ${\mathcal S}^{m-1}\left(\{z\in F:\ z\in\partial^*
E_w,\,w+z\notin\partial^*E\right\})=0$ for ${\mathcal S}^{k-m}$-a.e. $w\in
L$.
\end{lemma}

Thanks to this fact, we have that if $F\subset G\subset QX^*$ are two finite dimensional spaces,
then the relative essential boundary $\partial^*_FE$ of $E$ is contained, up to negligible sets, 
into the essential boundary $\partial^*_GE$ of $E$ relative to $G$, that is 
\[
{\mathcal S}^{\infty-1}_F(\partial^*_F E\setminus \partial^*_G E)=0.
\]

In \cite{AMP} there is the proof of the following fact.

\begin{proposition} \label{p2} Let ${\mathcal F}$ be a countable family of
finite-dimensional subspaces of $QX^*$ stable under finite unions.
For $F\in {\mathcal F}$, let $A_F\in{\mathcal B}(X)$ be such that
\begin{itemize}
\item[(i)] ${\mathcal S}^{\infty-1}_F(A_F\setminus A_G)=0$ whenever $F\subset G$;
\item[(ii)] $\sup_F{\mathcal S}^{\infty-1}_F(A_F)<\infty$.
\end{itemize}
Then $\displaystyle{\lim_F({\mathcal S}^{\infty-1}_F\res A_F)}$ exists, and it is
representable as $\displaystyle{(\lim_F{\mathcal S}^{\infty-1}_F)\res A}$ with
\[
A:=\bigcup_{F\in{\mathcal F}}\bigcap_{G\in{\mathcal F},\,
G\supset F}A_G\in{\mathcal B}(X).
\]
\end{proposition}

Such Proposition allows for the definition of the cylindrical essential boundary.

\begin{defn}[Cylindrical essential boundary] \label{defcyless}
Let ${\mathcal F}$ be a countable set of finite-dimensional
subspaces of $H$ stable under finite union, with
$\cup_{F\in{\mathcal F}}F$ dense in $H$. Then,
we define cylindrical essential boundary $\partial_{\mathcal F}^*E$
along ${\mathcal F}$ the set
\bigskip
\begin{equation*}
\partial_{\mathcal F}^*E:=\bigcup_{F\in{\mathcal
F}}\bigcap_{G\in{\mathcal F},\,G\supset F}\partial_G^*E.
\end{equation*}
\end{defn}

These definitions are used in \cite{hin09set} and with minor revisions in \cite{AMP}, to 
get a representation of the perimeter measure as follows.

\begin{theorem}\label{thmhino}
Let $E\in{\mathcal B}(X)$ be a set with finite $\gamma$-perimeter in
$X$, let ${\mathcal F}$ be as in Definition~\ref{defcyless} and let
$\partial^*_{\mathcal F}E$ be the corresponding cylindrical
essential boundary. Then
\begin{equation}\label{late2}
|D_\gamma\chi_E|(B)={\mathcal S}^{\infty-1}_{\mathcal
F}(B\cap\partial^*_{\mathcal F}E)\qquad\forall B\in{\mathcal B}(X).
\end{equation}
In particular, $\partial^*_{\mathcal F}E$ is uniquely determined by
\eqref{late2} up to ${\mathcal S}^{\infty-1}_{\mathcal F}$-negligible sets.
\end{theorem}

In \cite{AMP} also a weak rectifiability result of the cylindrical essential boundary is given; 
the term weak refers to the fact that rectifiability is done by using Sobolev functions instead 
of Lipschitz maps as in \eqref{countrectset}. This is not a minor difficulty, since in the 
infinite-dimensional setting no Lusin type properties are known; in particular, it is not known if
any Sobolev function coincides with a Lipschitz map in a set of positive measure.

First, we recall the notion of $H$-graph.

\begin{defn}[$H$-graph]
A set $\Gamma\subset X$ is called an $H$-graph if there exist
a unit vector $k\in QX^*$ and $u:D\subset {\rm Ker}(\pi_F)\to\R$
(here $F=\{sk,\ s\in \R\}$) such that 
\[
\Gamma=\{ y+u(y)k:\ y\in D\}.
\]
We say that $\Gamma$ is an entire Sobolev $H$-graph if moreover 
\[
D\in{\mathcal B}({\rm Ker}(\pi_F)),
\gamma^\perp\bigl({\rm Ker}(\pi_F)\setminus D\bigr)=0
\]
and $u\in W^{1,1}({\rm Ker}(\pi_F),\gamma^\perp)$.
\end{defn}

With this notion, in \cite{AMP} the following theorem is proved.

\begin{theorem}\label{thmBVHrect}
For any set $E\subset X$ with finite perimeter the measure $|D_\gamma\chi_E|$
is concentrated on a countable union of entire Sobolev $H$-graphs.
\end{theorem}

In \cite{AmbFig2}, the Ornstein-Uhlenbeck semigroup is used to define points 
of density $1/2$; their main result can be summarised in the following Theorem.

\begin{theorem}
Let $E\subset X$ be a set with finite perimeter; then 
\[
\lim_{t\downarrow 0} \int_X \left| T_t\chi_E -\frac{1}{2}\right|^2 d|D_\gamma\chi_E| =0;
\]
in particular, there exists a sequence $t_i\downarrow 0$ such that
\begin{equation}\label{esSerTt1}
\sum_{i} \int_X  \left| T_{t_i}\chi_E -\frac{1}{2}\right| d|D_\gamma\chi_E| <+\infty,
\end{equation}
which ensures that $T_{t_i}\chi_E\to \frac{1}{2}$ $|D_\gamma\chi_E|$-a.e. in $X$.
\end{theorem}

Thanks to the previous Theorem, a notion of points of density $\frac{1}{2}$ can be given. As 
explained in connection with the notion of essential boundary, the analogue \eqref{density1} 
of the finite dimensional procedure \eqref{defdensity} is not available in the present 
situation, hence it relies rather on an approach analogous to \eqref{heatdensity}. 

\begin{defn}[Points of density $1/2$]
Let $(t_i)_i$ be a sequence such that 
\begin{equation}\label{sqrtConv}
\sum_i \sqrt{t_i}<+\infty
\end{equation}
and \eqref{esSerTt1} holds. Then, we say that $x$ is a point of density $\frac{1}{2}$ for $E$ 
if it belongs to 
\begin{equation}\label{defDens12OU}
E^{1/2}:= \left\{
x\in X: \exists \lim_{i\to +\infty} T_{t_i}\chi_E(x)=\frac{1}{2}
\right\}.
\end{equation}
\end{defn}

The requirement in \eqref{sqrtConv} is rather natural, since for a set with finite perimeter it is possible to
prove (see \cite[Lemma 2.3]{AmbFig2}) that
\[
\int_X |T_t\chi_E-\chi_E| d\gamma \leq c_t P_\gamma(E),
\]
with
\[
c_t=\sqrt{\frac{2}{\pi}} \int_0^t \frac{e^{-s}}{\sqrt{1-e^{-2s}}}ds \sim 2\sqrt{\frac{t}{\pi}}.
\]

\begin{theorem}
Let $(t_i)_i$ be a sequence such that $\sum_i \sqrt{t_i}<+\infty$ and \eqref{esSerTt1} holds.
Then $|D_\gamma\chi_E|$ is concentrated on $E^{1/2}$ defined in \eqref{defDens12OU}; moreover 
$E^{1/2}$ has finite ${\mathcal S}^{\infty-1}$ measure and
\[
|D_\gamma\chi_E|={\mathcal S}^{\infty-1}\res E^{1/2}.
\]
\end{theorem}

It is worth noticing that the sequence $(t_i)_i$ depends of the set $E$ itself. In \cite{AmbFigRuna}it 
is also proved a part of the rectifiability result for the reduced boundary; with minor revision
of the definition of cylindrical essential boundary, it is possible to define a 
{\em cylindrical reduced boundary} by setting
\[
{\mathcal F}_F E=\{ x\in X: x=y+z: z\in {\mathcal F}(E_y) \subset F \},
\]
and
\begin{equation}\label{defRedBdryH}
{\mathcal F}_H E = \liminf_{F\in {\mathcal F}} {\mathcal F}_F E
=\bigcup_{F\in {\mathcal F}}\bigcap_{G\in{\mathcal F},G\supset F} {\mathcal F}_GE,
\end{equation}
where here ${\mathcal F}$ has two meanings, the first one to denote the reduced boundary, the
second one when writing $F\in {\mathcal F}$ is meant as a countable collection of finite dimensional sets
as in Proposition~\ref{p2}. The liminf of sets in \eqref{defRedBdryH} is also given in the sense of 
Proposition~\ref{p2}.

Given an element $h\in H$, the halfspace having $h$ as its ``inner normal"  is defined as
\[
S_h=\{ x\in X: \hat h(x)>0\}. 
\]
Notice that $S_h$ is a closed halfspace if $\hat{h}=R^*x^*$ for some $x^*\in X^*$; otherwise, it is easily 
seen by approximation that $\hat{h}$ is linear on a subspace of $X$ of full measure, hence the above 
definition does make sense. Since the convergence of sequences $h_n\in H$ to $h\in H$ in the norm of 
$H$ implies the convergence of $S_{h_n}$ to $S_h$ in the sence of convergence of characteristic functions 
in $L^1(X,\gamma)$, then, denoting by
\[
E_{x,t}:=\frac{E-e^{-t}x}{\sqrt{1-e^{-2t}}},
\]
the following result holds true. We notice that the idea underlying the following result 
is the last line in \eqref{OUsgrpinRd}, which cannot be used directly in the infinite-dimensional 
framework. 

\begin{theorem}[Ambrosio, Figalli, Runa~\cite{AmbFigRuna}]
Let $E\subset X$ be a set with finite perimeter in $X$, $x\in {\mathcal F}_HE$ and $S(x)=S_{\nu_E(x)}$
where $\nu_E$ is defined by the polar decomposition $D_\gamma\chi_E=\nu_E |D_\gamma\chi_E|$;
then
\[
\lim_{t\downarrow 0} \int_X \int_X
\left|\chi_E(e^{-t}x+\sqrt{1-e^{-2t}}y)-\chi_{S(x)}(y)\right| d\gamma(y) d|D_\gamma\chi_E|(x)=0.
\] 
\end{theorem}

In other terms, the previous results can be restated by saying that
\[
\lim_{t\downarrow 0} \int_X \|\chi_{E_{x,t}}-\chi_{S(x)}\|_{L^1(X,\gamma)} d|D_\gamma\chi_E|(x)=0,
\]
that is, $E_{x,t}$ coverge to $S(x)$ in $L^1(X,\gamma)$, for $|D_\gamma\chi_E|$-a.e. $x\in X$.
This result is in some sense the Wiener space formulation of \eqref{apptgredbdry}. 

\subsection{Examples of sets with finite perimeter}

We now provide some examples of sets with finite perimeter; in some cases the essential and 
reduced boundary are directly identifiable, in some other they are indicated as candidates, but 
a proof is not available so far.

\subsubsection{Cylindrical sets.}
Let ${\mathcal F}$ be as in Definition \ref{defcyless}. The easiest way to construct examples of 
sets with finite perimeter is to use the decomposition $X=X_m\oplus {\rm Ker}(\pi_F)$; if $B\subset F$ 
is a set with $\chi_B\in BV(X_m,\gamma_F)$, then $E=\pi_F^{-1}(B)$ has finite perimeter in $X$ with 
\[
P_\gamma(E,X)=P_{\gamma_F}(B,X_m).
\]
If $F\in {\mathcal F}$, then 
\[
\partial^*_{\mathcal F}E =\partial^*_F E=\partial^* B,\qquad
{\mathcal F}_H E={\mathcal F} B,
\]
otherwise the previous equality holds up to $|D_\gamma\chi_E|$-negligible sets.

\subsubsection{Level sets of Lipschitz maps: comparison with the Airault-Malliavin surface measure.} 
By coarea formula \eqref{coareaform}, almost every level set of a $BV$ function has finite perimeter; 
in particular, we can use almost every level set of Sobolev or Lipschitz functions.
To prove that {\em every} level set, under some regularity assumption on the function, has finite perimeter
is quite delicate in this framework. In \cite{AirMal}, Airault and Malliavin constructed a surface measure 
on boundaries of regular level sets. More precisely, they considered functions $f$ belonging to
\[
W^\infty(X,\gamma)=\bigcap_{p>1, k\in \N} W^{k,p}(X,\gamma),
\]
where $W^{k,p}(X,\gamma)$ is the Sobolev space of order $k$ with $p$-integrability, 
such that
\[
\frac{1}{|\nabla_H f|_H} \in \bigcap_{p\geq 1} L^p (X,\gamma);
\]
what they proved is that the image measure $f_\# \gamma$ defined on ${\mathcal B}(\R)$ by
\[
f_\# \gamma(I)=\gamma(f^{-1}(I))
\]
has smooth density $\rho$ with respect to the Lebesgue measure and that, for each $t$ such that
$\rho(t)>0$, there exists a Radon measure $\sigma_t$ supported on $f^{-1}(t)$ such that
\[
\int_{\{f<t\}} \diver_H \Phi\, d\gamma =\int_{\{f=t\}} \frac{[\Phi,\nabla_H f]_H}{|\nabla_H f|_H}\, d\sigma_t.
\]
The measure $\sigma_t$ is constructed in terms of the Minkowski content as explained in 
Subsection \ref{isoper_ineq}. In \cite{CasLunMirNov}, it is proved that, under the additional 
technical assumption that $f$ is continuous, the set $\{f<t\}$ has finite
perimeter whenever $\rho(t)>0$ with the identity
\[
P_\gamma(\{f<t\})=\sigma_t(\{f=t\})=\int_{\{f<t\}} \diver_H \nu_H \, d\gamma,
\]
where $\nu_H=\nabla_H f/|\nabla_H f|_H$. The set $\{f=t\}$ is expected to be the essential 
boundary of $\{f<t\}$, whereas the points in the reduced boundary are expected to be those
$x$ where $\nabla_H f(x)\neq 0$. 

\subsubsection{Balls and convex sets.}
If we fix a point $x_0\in X$, the map
\[
f(x)=\|x-x_0\|_X
\]
is Lipschitz and then the sets
\[
E_t=\{ f<t\}=B_t(x_0)
\]
have finite perimeter for almost every $t>0$. The proof that {\em every} ball has finite perimeter is 
contained in \cite{CasLunMirNov}; if $X$ is a Hilbert space, then the function $f(x)^2$ is continuous 
and satisfies all the condition imposed by Airault and Malliavin and then all balls in Hilbert spaces 
have finite perimeter. In addition, the normal vector in this case is given by
\[
\nu(x) =\frac{Q(x-x_0)}{|Q(x-x_0)|_H}
\]
(where $Q$ is the covariance operator, $\gamma=\NN(0,Q)$) and the function
\[
g(t)= P_\gamma(B_t(x_0))
\]
is continuous in $[0,+\infty)$ with
\[
\lim_{t\to 0}P_\gamma(B_t(x_0))=\lim_{t\to +\infty} P_\gamma(B_t(x_0))=0.
\]
It is also possible to prove that there exist $t_1<t_2$ such that $g$ is increasing in $[0,t_1]$
and decreasing in $[t_2,+\infty)$.

The proof that any ball in an infinite-dimensional Banach space has finite perimeter
is less explicit and is based on a Brunn-Minkowski argument stating that for every Borel sets $A,B\subset X$, 
\[
\gamma(\lambda A+(1-\lambda )B)\geq \gamma(A)^\lambda \gamma(B)^{1-\lambda},\qquad
\lambda\in [0,1].
\] 
In \cite{CasLunMirNov} it is proved that if $C$ is an open convex set, then $\gamma(\partial C)=0$ and $C$ 
has finite perimeter.
In this case, it is easily seen that $\partial^*_{\mathcal F}C\subset \partial C$ and 
\begin{equation}\label{brdConv}
|D_\gamma\chi_C|(\partial C\setminus \partial^*_{\mathcal F}C)=0;
\end{equation}
indeed if $x\in C$ or $x\in X\setminus \overline C$, then for any $F\leq H$, if we write 
\[
x=y+z_x,\qquad y\in {\rm Ker}(\pi_F), z_x\in X_m
\]
then $z_x$ is an interior point either of $C_y$ or of $X_m\setminus \overline C_y$, so
$\partial^*_{\mathcal F}C\subset \partial C$. Property \eqref{brdConv} follows by the representation
of the perimeter measure \eqref{late2}. The characterization of the reduced cylindrical boundary is 
less clear.

The assumption that $C$ is open is essential; indeed, it is also shown that, in the Hilbert space case, 
there exists a convex set with infinite perimeter. Such a set is constructed by fixing a sequence $r_i$ 
such that
\[
\sqrt{\frac{2}{\pi}}\frac{e^{-\frac{r_i^2}{2}}}{r_i}=\frac{1}{(i+1)(\log (i+1))^\frac{3}{2}},
\]
defining
\[
C_m=\pi_F^{-1}(Q_m), \quad Q_m=\prod_{i=1}^m [-r_i,r_i]
\]
and letting $m\to +\infty$.

\subsubsection{An example in the classical Wiener space.}

In \cite{HinUch08} an example of a set with finite perimeter in the classical Wiener space is given, 
using the reflecting Brownian motion. The setting is given by a {\em pinned path space}, that is 
\[
X=\{ \omega\in C([0,1],\R^d): \omega(0)=a, \omega(1)=b\}
\]
endowed with the pinned Wiener measure  $\PP_{a,b}$ defined in the same spirit as 
\eqref{Wienermeas} by
\[
\PP_{a,b}(C) = \frac{1}{G(1,a,b)}\int_{B_1\times\ldots \times B_m} 
\prod_{j=1}^{m+1} G(t_j-t_{j-1},x_{j-1},x_j)\, dx_1 \ldots dx_m,
\]
where $B_j\in \B(\R^d)$, $j=1,\ldots, m$, $0=t_0<t_1<\ldots <t_m < t_{m+1}=1$,
$x_0=a$ and $x_{m+1}=b$, 
\[
C=\{\omega\in X:\ \omega(t_j)\in B_j,\ j=1,\ldots, m\}.
\]
In such space, if $\Omega\subset \R^d$ is an open set containing the two points $a$ and $b$, define the set
\[
E^\Omega=\{ \omega\in X: \omega(t)\in \Omega\ \forall t\in [0,1]\}.
\]
Then $E^\Omega$ has finite perimeter in $X$ under the assumption that $\Omega$ has positive reach, that is
an uniform exterior ball condition: there exists $\delta>0$ such that for every $y\in \partial \Omega$ 
there is $z\in \R^d\setminus \Omega$ such that $B_\delta(z)\cap \overline \Omega=\{ y\}$. The proof of 
this fact is done constructing a sequence of Lipshitz functions $\rho_n$ converging
to $\chi_{\overline E^\Omega}$ in $L^1(X,\PP_{a,b})$ and such that
\[
\int_X |\nabla_H\rho_n|_H d\PP_{a,b} \leq n \PP_{a,b}\left(
\left\{ \omega\in X: 0\leq \inf_{t\in [0,1]} q(\omega(t))\leq \frac{1}{n} \right\}
\right);
\]
the sequence is defined in terms of the signed distance function
\[
q(x)=\inf_{\inf y\in \R^d\setminus \Omega} |x-y| -\inf_{y\in \Omega} |y-x|
\]
as 
\[
\rho_n(\omega) =f_n(F(\omega)),\quad F(\omega) =\inf_{t\in [0,1]} q(\omega(t)),
\]
where $f_n$ is defined as
\[
f_n(s)=\min \{ \max\{0,ns\},1\}.
\]
The keypoint in the proof where the positive reach condition is used is in estimating
\[
\PP_{a,b}\left(\Big\{ \omega\in X: 0\leq \inf_{t\in [0,1]} q(\omega(t))\leq r\Big\}\right)\leq cr,
\]
since from that it comes that
\[
\int_{X_{a,b}} |\nabla_H \rho_n |_H d\PP_{a,b} \leq c.
\]
In this case,  Hino-Uchida prove also that the perimeter measure concentrates on the set
\[
\partial'E^\Omega =\left\{\omega\in X: \omega(t)\in \overline \Omega \mbox{ and }
\exists \mbox{ an unique t }\in [0,1] \mbox{ s.t. } \omega(t)\in \partial \Omega
\right\}.
\]
The definition of the previous set has a meaning very close to the set of points of density $1/2$ for $E^\Omega$.
Finally, it is worth noticing that the proof given by Hino and Uchida of the fact that $E^\Omega$ ha finite 
perimeter is close to the proof that a (sufficiently regular) set in the Euclidean setting has finite 
Minkowski content. 

\section{Convex functionals on $BV$}\label{secfunctionals}

Following \cite{ChaGolNov}, we now consider integral functionals on $BV(X,\gamma)$ of the form
\[
u\mapsto \int_X F(D_\gamma u) 
\]
where $F:H \to \R\cup\{+\infty\}$ is a convex lower semicontinuous function.
As $D_\gamma u$ is in general a measure, we have to give a precise meaning to the above expression.

Given a convex function $F:H \to \R\cup \{+\infty\}$  we denote by $F^*$ its convex conjugate, defined as
\[
F^*(\Phi) :=\sup\left\{ [\Phi,h]_H -F(h):  h \in H\right\}, \qquad \Phi \in H,
\]
and by $F^\infty$ its recession function defined  as
\[
F^\infty(h):=\lim_{t \to +\infty} \frac{F(th)}{t} \qquad h \in H.
\]
We shall consider functions $F: H\to \R\cup\{+\infty\}$ satisfying the following assumption:
\begin{itemize}
\item[(A)] $F$ is a proper (i.e., not identically $+\infty$), lower semi-continuous, convex function on $H$.
\end{itemize}
\noindent Notice that a convex function $F$ with $p\ge 1$ growth, i.e., such that there are positive constants 
$\alpha_1$, $\beta_1$, $\alpha_2$, $\beta_2$ such that
\begin{equation}\label{condacca}
\alpha_1 \normH{h}^p-\beta_1\le F(h)\le \alpha_2 \normH{h}^p +\beta_2 \qquad \forall h \in H,
\end{equation}
satisfies automatically assumption (A).

Given a function $F$ satisfying (A) and $u\in \Ldeu$, we define the functional 
\begin{equation}\label{repsupF}
\int_X F(D_\gga u):=\sup\Bigl\{ \int_X -u \diver_H \Phi - F^*(\Phi) \ d\gga,\ 
\Phi \in \FCb(X,H)\Bigr\}
\end{equation}
which is lower semicontinuous in $\Ldeu$. 
Similarly, for $\mu \in \cM(X,H)$ we set
\[\int_X F(\mu) := \sup\Bigl\{ \int_X [\Phi, d\mu]_H -\int_X F^*(\Phi)d\gga,\ 
\Phi \in \FCb(X,H)\Bigr\}.
\]
The following result has been proved in \cite[Theorem 3.2]{ChaGolNov}.

\begin{theorem}\label{corDemTem}
Let $F : H \to \R\cup \{+\infty\}$ satisfy (A) and let $\mu \in \cM(X,H)$, then  
\begin{equation*}
\int_X F(\mu)= \int_X F(\mu^a) d\gga+\int_X F^{\infty}\left(\frac{d \mu^s}{d|\mu^s|} \right)d|\mu^s|
\end{equation*}
where $\mu=\mu^a \gga +\mu^s$ is the Radon-Nikodym decomposition of $\mu$ w.r.t. $\gamma$.
\end{theorem}

{}From Theorem \ref{corDemTem} we obtain a representation result for the functional in \eqref{repsupF}.

\begin{theorem}\label{thF}
Let $F : H \to \R\cup \{+\infty\}$ satisfy (A), then
\begin{equation}
\int_X F(D_\gga u)= 
\int_X F(\nabla_H u) d\gga+ \int_X F^\infty\left(\frac{d D^s_\gga u}{d|D^s_\gga u|}\right) d|D^s_\gga u| 
\end{equation}
for all $u \in BV(X,\gga)$, where $D_\gga u =\nabla_H u \gga +D^s_\gga u$ is the Radon-Nikodym decomposition 
of $D_\gga u$. 
\end{theorem}

A natural question is whether the functional  in \eqref{repsupF} concides with the relaxation in $\Ldeu$
of its restrictions to more regular functions. 
The following result has been proved in \cite[Proposition 3.4]{ChaGolNov}.

\begin{theorem}\label{relax}
Let $F : H \to \R\cup \{+\infty\}$ satisfy (A), then the functional $\int_X F(D_\gga u)$ is 
the relaxation in $\Ldeu$ of the functional defined as $\int_X F(\nabla_H u) d\gga$ for 
$u \in W^{1,1}(X,\gga)$, and $+\infty$ for $u \not\in W^{1,1}(X,\gga)$.

If $F$ has $p\ge 1$ growth in the sense of \eqref{condacca},  
then the same relaxation result holds with the space $W^{1,1}(X,\gga)$ replaced by $\FCb(X)$.
\end{theorem}

Condition \eqref{condacca} in the above statement is technical, and we expect that it is not 
necessary to obtain the relaxation result in $\FCb(X)$.

\subsection{Convexity of minimisers}

The Direct Method of the Calculus of Variations is a well-known method to prove existence of 
minimisers of variational problems.
The two conditions a functional has to satisfy in order to apply the method are the lower 
semicontinuity with respect to a given topology,
and the compactness of a nonempty sublevel set in the same topology. 

We now consider convex functionals of the form
\begin{equation}\label{eqcon}
\int_X F(D_\gga u) +\frac{1}{2}\int_X (u-g)^2 d\gga.
\end{equation}
where $F: H\to \R\cup\{+\infty\}$ satisfies (A) and $g \in \Ldeu$ is a convex function. 

Notice that the functional in \eqref{eqcon} is convex on $\Ldeu$,
hence it is also weakly lower semicontinuous. Moreover,
its sublevel sets are (relatively) compact in the weak 
topology of $\Ldeu$. By the Direct Method we then obtain the 
following existence result. The existence of a minimiser follows by the Direct Method 
of the Calculus of Variations, while the uniqueness follows from the strict convexity 
of the functional, due to the second term in \eqref{eqcon}.

\begin{proposition}\label{proex}
There exists a unique minimiser $\bar u\in \Ldeu$ of the functional \eqref{eqcon}.
\end{proposition}

We state a convexity result for minimisers of \eqref{eqcon} which has been proved in 
\cite[Theorem 5.1]{ChaGolNov}.

\begin{theorem}\label{convlinear}
The minimiser $\bar u$ of \eqref{eqcon} is convex.
\end{theorem}

{}From Theorem \ref{convlinear} and the theory of maximal monotone operators (see \cite{brezis}), 
one can easily get the following result:

\begin{theorem}
Let $u_0\in \Ldeu$ be a convex initial datum.
Then the solution $u(t)$ of the $\Ldeu$-gradient flow of  $\int_X F(D_\gga u)$ with initial 
condition $u(0)=u_0$ is convex for every $t>0$.
\end{theorem}

Notice that, by taking $F(h)=|h|^p$ with $p\ge 1$, 
Theorem \ref{convlinear} applies to the functional
\begin{equation}\label{eqmin}
\int_X \normH{D_\gga u}^p +\frac{1}{2} \int_X (u-g)^2 d\gga.
\end{equation}
Recalling the coarea formula \eqref{coareaform}, when $p=1$ the functional \eqref{eqmin} can be written as
\[
\int_X \normH{D_\gga u}^p +\frac{1}{2} \int_X (u-g)^2 d\gga = 
\int_\R \Bigl( P_\gga(\{u>t\}) - \int_{\{u>t\}}\!\!\!(g-t)d\gga\Bigr) dt.
\]
It then follows (see \cite{CasMirNov1,ChaGolNov}) that the level set $\{\bar u>t\}$ of the minimiser $\bar u$ 
minimises the geometric problem
\begin{equation}\label{fungeo}
P_\gga(E) - \int_{E}(g-t)\,d\gga
\end{equation}
among the subsets $E\subset X$ of finite perimeter, for all $t\in\R$. 
Then, from Theorem \ref{convlinear} one can derive a convexity result 
for minimisers to \eqref{fungeo} (see \cite[Corollary 5.7]{ChaGolNov}).

\begin{theorem}
Let $g\in\Ldeu$ be a convex function, and consider the functional 
\begin{equation}\label{fungeo2}
F_g(E) = P_\gga(E) - \int_{E}g\,d\gga.
\end{equation}
Then, two situations can occur:
\begin{itemize}
\item If $\min F_g < 0$, there exists a unique nonempty minimiser of $F_g$, which is convex.
\item If $\min F_g = 0$, there exists at most one nonempty minimiser of $F_g$, which is then convex.
\end{itemize}
\end{theorem}

\subsection{Relaxation of the perimeter in the weak topology}\label{secperimeter}

In view of the previous discussion, a natural problem which arises is the 
classification of the weakly lower semicontinuous functionals on $\Ldeu$. 

While convex functionals are lower semicontinuous with respect to both the weak and the strong topology,
the perimeter functional 
\[
F(u):=
\left\{\begin{array}{ll}
\Pgamma{E} \qquad& \textrm{if } u=\chi_E\\[8pt]
+\infty \qquad & \textrm{otherwise}
\end{array}\right.
\]
is not weakly lower semicontinuous, as one can easily check by taking the sequence of 
halfspaces ${E_n}=\{ \sprod{x}{x_n^*}<0\}$, where $x_n^*$ is a sequence in $X^*$ such that 
$h_n = Q x_n^*$ is an orthonormal basis of $H$. Indeed, the characteristic functions of these sets weakly 
converge to the constant function $1/2$, which is not a characteristic function,
while the perimeter of $E_n$ is constantly equal to $1/\sqrt{2\pi}$.

In \cite{GolNov} the authors computed the relaxation ${\overline F}$ of $F$ with respect to the 
weak $\Ldeu$-topology, showing that 
\[
{\overline F}(u)=\left\{\begin{array}{ll}
\displaystyle \int_X \sqrt{\cU^2(u)+|D_\gga u|^2} \qquad& 
\textrm{if } u \in BV(X,\gga) \textrm{ and } |u|\le 1 
\\[8pt]
+\infty \qquad & \textrm{otherwise}
\end{array}\right.
\]
where 
\begin{equation}\label{intrep}
\int_X \sqrt{\cU^2(u)+|D_\gga u|^2}=\int_X \sqrt{\cU^2(u)+\normH{\nabla_H u}^2}d \gga + |D^s_\gga u|(X)
\end{equation} 
with $D_\gga u= \nabla_H u \,d \gga + D^s_\gga u$ as in Theorem \ref{thF}. 
Observe that the functional ${\overline F}$ already appears in the seminal works by Bakry and Ledoux 
\cite{bakryledoux}  and Bobkov \cite{bobkov}, in the context of log-Sobolev inequalities. 
See also \cite[Remark 4.3]{AMMP} where it appears in a setting closer to ours.

There is also a representation formula for ${\overline F}$, which is reminiscent of the definition 
of total variation:
\begin{align*} 
{\overline F}(u)= \sup\Bigl\{&\int_X (u \diver_H \Phi + \cU(u) \xi ) d\gga \,:\ 
\Phi\in \FCb(X,H),\ 
\\
&\xi \in \FCb(X),\ \normH{\Phi(x)}^2+|\xi(x)|^2\le 1 \; \forall x \in X \Bigr\},
\end{align*}
for all $u \in BV(X,\gga)$, with $|u|\le 1$.

\section{Open problems} \label{openproblems}

We collect some open problems whose solution, in our opinion, would provide important 
information on the whole subject and would allow for a wide range of applications. 

The first problems that should be solved and would have a great influence in the further 
developments concern the structure theory of reduced boundaries and general $BV$ functions. 
For instance it would be important to check whether the well-known Euclidean decomposition result  
holds in Wiener spaces, i.e., whether the equality $X=E^1\cup E^0\cup E^{1/2}$ is true
(up to negligible sets). Moreover, as  
we have seen, a {\em pointwise} characterization of reduced boundary like that in 
\eqref{redbdry} is missing, as well as suitable notions of {\em one-sided approximate limits}, 
see \eqref{halflimits}. In this respect, the Orstein-Uhlenbeck semigroup will come into play, 
but making density computations independent of the sequence $(t_i)$, see \eqref{defDens12OU}, 
would certainly be useful, in connection with the coarea formula. Still on the side of the 
structure theory, it is important to improve the weak rectifiability Theorem \ref{thmBVHrect},
possibly getting Lipschitz rectifiability. All these problems are of course connected to the 
general problem of the traces of $BV$ functions. Beside other instances, such as boundary 
value problems, closer to the arguments presented here are applications of the structure 
theory and fine properties to integral functionals. Indeed, it would be interesting to extend 
the results presented in Section \ref{secfunctionals} to integrands depending on $u$, see 
\cite[Section 5.5]{afp} for the classical case. In this connection, it would be important to 
perform a deeper analysis of the singular part of the gradient, possibly distinguishing between 
the {\em jump part} and the {\em Cantor part}, and defining the one-sided approximate limits. 
This could probably give a representation formula more precise than \eqref{intrep}.  
Finally, one could try to provide a complete characterization of weakly lower semicontinuous 
integral functionals with integrands of linear growth.



\begin{thebibliography}{99}

\bibitem{AirMal}
H. Airault, P. Malliavin, 
Int\'egration g\'eom\'etrique sur l'espace de Wiener, 
Bull. Sci. Math. {\bf 112(1)} (1988), 3-52.

\bibitem{ADP}
L. Ambrosio, G. Da Prato, D. Pallara, 
$BV$ functions in a Hilbert space with respect to a Gaussian measure, 
Rend. Acc. Lincei {\bf 21}, (2010), 405-414. 

\bibitem{ADPGP}
L. Ambrosio, G. Da Prato, B. Goldys, D. Pallara, 
Bounded variation with respect  to a log-concave measure, 
Comm. P.D.E., to appear.

\bibitem{AmbFig1}
L. Ambrosio, A. Figalli,
On flows associated to Sobolev vector fields in Wiener spaces: an approach \`a la Di Perna-Lions,
J. Funct Anal.  {\bf 256}, (2009), 179-214.

\bibitem{AmbFig2}
L. Ambrosio, A. Figalli,
Surface measure and convergence of the Ornstein-Uhlenbeck semigroup in Wiener spaces, 
Ann. Fac. Sci. Toulouse Math., {\bf 20}, (2011) 407-438.

\bibitem{AmbFigRuna}
L. Ambrosio, A. Figalli, E. Runa, 
On sets of finite perimeter in Wiener spaces: reduced boundary and convergence to halfspaces,
Preprint (2012).

\bibitem{afp} 
L. Ambrosio, N. Fusco, D. Pallara, 
Functions of bounded variation and free discontinuity problems, 
Oxford Mathematical Monographs, 2000.

\bibitem{AMMPPhysicad}
L.  Ambrosio, S. Maniglia, M.  Miranda Jr, D.  Pallara, 
Towards a theory of $BV$ functions in abstract Wiener spaces, 
In: ``Evolution Equations: a special issue of Physica D'', 
{\bf 239} (2010), 1458-1469. 

\bibitem{AMMP}
L.  Ambrosio, S. Maniglia, M.  Miranda Jr, D.  Pallara, 
$BV$ functions in abstract Wiener spaces, J. Funct. Anal., 
{\bf 258} (2010), 785-813.

\bibitem{AMP}
L. Ambrosio, M. Miranda Jr, D. Pallara, 
Sets with finite perimeter in Wiener spaces, perimeter measure and boundary rectifiability, 
Discrete Contin. Dyn. Syst., {\bf 28} (2010), 591-606.

\bibitem{bakryledoux} 
D.~Bakry, M.~Ledoux, 
L\'evy-Gromov's isoperimetric inequality for an infinite dimensional diffusion generator, 
Invent. Math., {\bf 123} (1996), 259-281.

\bibitem{bobkov} 
S.G.~Bobkov, 
An isoperimetric inequality on the discrete cube, and an elementary proof of the isoperimetric 
inequality in Gauss space, 
Ann. Probab., {\bf 25} (1997), 206-214. 

\bibitem{boga} 
V. I. Bogachev, 
Gaussian Measures, Mathematical Surveys and Monographs {\bf 62}, 
American Mathematical Society, 1998.

\bibitem{bogamt}
V. I. Bogachev, 
Measure Theory, Springer, 2007.

\bibitem{bogadm}
V. I. Bogachev,
Differentiable measures and the Malliavin calculus,
Mathematical Surveys and Monographs {\bf 164},
American Mathematical Society, 2010.  

\bibitem{brezis}
H.~Br\'ezis, 
Op\'erateurs maximaux monotones et semi-groupes de contractions dans les espaces de Hilbert, 
North Holland, 1973.

\bibitem{CarKer01OnT}
E.A.~Carlen and C.~Kerce,
On the cases of equality in Bobkov's inequality and Gaussian rearrangement, 
Calc. Var., 13 (2001), 1-18.

\bibitem{CasLunMirNov}
V. Caselles, A. Lunardi, M. Miranda jr, M. Novaga,
Perimeter of sublevel sets in infinite dimensional spaces,
Adv. Calc. Var., {\bf 5} (2012), 59-76.

\bibitem{CasMirNov1}
V. Caselles, M. Miranda jr, M. Novaga, 
Total Variation and Cheeger sets in Gauss space, 
J. Funct. Anal., {\bf 259} (2010), 1491-1516.

\bibitem{CNV}
A. Cesaroni, M. Novaga, E. Valdinoci,
A symmetry result for the Ornstein-Uhlenbeck operator,
Discrete Contin. Dyn. Syst. A, to appear.

\bibitem{ChaGolNov}
A. Chambolle, M. Goldman, M. Novaga, 
Representation, relaxation and convexity for variational problems in Wiener spaces, 
J. Math. Pures Appl., to appear.

\bibitem{DaPGolZab}
G. Da Prato, B. Goldys, J. Zabczyk, 
Ornstein-Uhlenbeck semigroups in open sets of Hilbert spaces, 
C. R. Math. Acad. Sci. Paris {\bf 325} (1997) 433-438.

\bibitem{DaPLun}
G. Da Prato, A. Lunardi,
On the Dirichlet semigroup for Ornstein-Uhlenbeck operators in subsets of Hilbert spaces, 
J. Funct. Anal. {\bf 259} (2010) 2642-2672. 

\bibitem{DeG1}
E. De~Giorgi, 
Su una teoria generale della misura $(r-1)$-dimensionale in uno spazio ad $r$ dimensioni, 
Ann. Mat. Pura Appl. (4), {\bf 36} (1954), 191-213.  Also in: Ennio De~Giorgi: {\em Selected Papers}, 
(L. Ambrosio, G. Dal Maso, M. Forti, M. Miranda, S. Spagnolo eds.), Springer, 2006, 79-99. 
English translation, {\em Ibid.}, 58-78.

\bibitem{DeG2}
E. De~Giorgi, 
Su alcune generalizzazioni della nozione di
perimetro, in: Equazioni differenziali e calcolo delle
variazioni (Pisa, 1992), G. Buttazzo, A. Marino, M.V.K. Murthy
eds, Quaderni U.M.I. {\bf 39}, Pitagora, 1995, 237-250.

\bibitem{Ehr83Sym}
A.~Ehrhard, 
Sym\'etrisation dans l'espace de {G}auss, 
Math. Scand., {\bf 53} (1983), 281-301.

\bibitem{evgar} 
L.C. Evans, R.F. Gariepy, 
{\em Lecture notes on measure theory and fine properties of functions,}
CRC Press, 1992.

\bibitem{federer} 
H. Federer, 
Geometric measure theory, Springer, 1969.

\bibitem{feydelpra}
D. Feyel, A. de la Pradelle, 
Hausdorff measures on the Wiener space, 
Potential Anal., {\bf 1} (1992), 177-189.

\bibitem{friedman}
A. Friedman, 
Stochastic differential equations and applications, 
Academic Press, 1975/76 and Dover, 2006.

\bibitem{fuk99}
M. Fukushima, 
On semimartingale characterization of functionals of symmetric Markov processes, 
Electron. J. Probab., {\bf 4} (1999), 1-32.

\bibitem{fuk2000_1}
M. Fukushima, 
$BV$ functions and distorted Ornstein-Uhlenbeck processes over the abstract Wiener space, 
J. Funct. Anal., {\bf 174} (2000), 227-249.

\bibitem{fuk2000_2}
M. Fukushima, M. Hino, 
On the space of BV functions and a Related Stochastic Calculus in Infinite Dimensions. 
J. Funct. Anal., {\bf 183} (2001), 245-268.

\bibitem{Fot94} 
M. Fukushima, Y. Oshima, M. Takeda, 
{\em Dirichlet forms and symmetric Markov processes,}
de Gruyter Studies in Mathematics {\bf 19}, de Gruyter, 1994.

\bibitem{GolNov}
M. Goldman, M. Novaga, 
Approximation and relaxation of perimeter in the Wiener space,
Annales IHP - Analyse Nonlineaire, {\bf 29} (2012), 525-544.

\bibitem{Hariya}
Y. Hariya,
Integration by parts formulae for Wiener measures restricted to subsets in $\R^d$, 
J. Funct. Anal. {\bf 239} (2006) 594-610.

\bibitem{hin04int}
M. Hino, 
Integral representation of linear functionals on vector lattices and its application to $BV$ functions 
on Wiener space, in: Stochastic analysis and related topics in Kyoto in honour of Kiyoshi It\^o, 
Advanced Studies in Pure Mathematics, {\bf 41} (2004), 121-140.

\bibitem{hin09set}
M. Hino, 
Sets of finite perimeter and the Hausdorff-Gauss measure on the Wiener space, 
J. Funct. Anal., {\bf 258} (2010), 1656-1681.

\bibitem{HinUch08} 
M.Hino and H.Uchida,
Reflecting Ornstein-Uhlenbeck processes on pinned path spaces,
Proceedings of RIMS Workshop on Stochastic Analysis and Applications, 111-128, 
RIMS Kokyuroku Bessatsu, B6, Kyoto, 2008. 

\bibitem{ledoux} 
M. Ledoux, 
Isoperimetry and Gaussian analysis, in: Lectures on Probability Theory and Statistics, 
Saint Flour, 1994, Lecture Notes in Mathematics {\bf 1648}, Springer, 1996, 165-294.

\bibitem{MaRockner}
Z. M. Ma and M. R\"ockner,
{\em Introduction to the theory of (non-symmetric) Dirichlet forms,}
Springer, 1992.

\bibitem{maggi}
F. Maggi, 
{\em Sets of finite perimeter and geometric variational problems: an introduction to geometric measure theory,} 
Cambridge Studies in Advanced Mathematics {\bf 135}, Cambridge University Press, 2012. 

\bibitem{mall76}
P. Malliavin,
Stochastic calculus of variation and hypoelliptic operators,
Diff. Eq. Res. Inst. Math. Sci., Kyoto Univ., Kyoto 1976,
Wiley, New York, (1978), 195-263.

\bibitem{malliavin}
P. Malliavin, 
Stochastic analysis, Grundlehren der Mathematischen Wissenschaften {\bf 313},
Springer, 1997.

\bibitem{preiss} 
D. Preiss,
Differentiability of Lipschitz functions in Banach spaces, 
J. Funct. Anal. {\bf 91}, (1990), 312-345.

\bibitem{Preiss} 
D. Preiss,  
Gaussian measures and the density theorem.
Comment. Math. Univ. Carolin., {\bf 22} (1981), 181-193.

\bibitem{sudtsi}
V.N. Sudakov and B.S. Tsirel'son, 
Extremal properties of half-spaces for spherically invariant measures. 
Problems in the theory of probability distributions, II, 
{\em Zap. Nau\v cn. Sem. Leningrad. Otdel. Mat. Inst. Steklov. (LOMI).},
{\bf 41}, (1974), 14-24.

\bibitem{trevisan}
D. Trevisan, 
$BV$-capacities on Wiener spaces and regularity of the maximum of the Wiener process, 
preprint. 

\bibitem{VTC} 
N. N. Vakhania, V. I. Tarieladze, S. A. Chobnyan, 
{\em Probability distribution in Banach spaces,}
Kluwer, 1987.

\bibitem{Volpert}
A. I. Vol'pert,
The space $BV$ and quasilinear equations,
Math. USSR Sbornik {\bf 2} (1967), 225-267. 

\bibitem{Wentzell}
A. D. Wentzell, 
{\em A Course in the Theory of Stochastic Processes,} 
McGraw-Hill, New York, 1981. 

\bibitem{Zambotti}
L. Zambotti, 
Integration by parts formulae on convex sets of paths and applications to SPDEs with reflection,
Probab. Theory Related Fields {\bf 123} (2002) 579-600.

\end{thebibliography}
\end{document}